\documentclass{elsarticle}
\makeatletter
\def\ps@pprintTitle{%
	\let\@oddhead\@empty
	\let\@evenhead\@empty
	\def\@oddfoot{}%
	\let\@evenfoot\@oddfoot}
\makeatother

\usepackage[utf8]{inputenc}
\usepackage{fontenc}[t1]
\usepackage[hidelinks]{hyperref}
\usepackage{mathptmx} 
\usepackage{setspace}
\usepackage{amssymb}
\usepackage{graphicx}
\usepackage{subfig}
\usepackage{caption}
\usepackage{amsfonts}
\usepackage{hyperref}
\usepackage{amsmath}
\usepackage{amsthm} 
\theoremstyle{plain} 
\usepackage{comment}
\usepackage{natbib}

\newcommand{\field}[1]{\mathbb{#1}}

\newcommand {\N}{ {\field{N}} }

\usepackage{color}
\usepackage{colordvi}

\newtheorem{prop}{Proposition}[section]

\newtheorem{cor}[prop]{Corollary}

\newcommand{\komment}[1]{}

\usepackage[left]{lineno}

\begin{document}

\begin{frontmatter}
	
	\title{Parameter estimation for contact tracing in graph-based models}
	
	\author[TUMGarching]{Augustine Okolie\corref{mycorrespondingauthor}} 
	\author[TUMGarching,ICB]{Johannes M\"uller} 
	\author[UC]{Mirjam Kretzschmar} 
	
	\address[TUMGarching]{Center for Mathematical Sciences, Technische Universit\"at M\"unchen, 85748 Garching, Germany}
	\address[ICB]{Institute for Computational Biology, Helmholtz Center Munich, 85764 Neuherberg, Germany}
	\address[UC]{University Medical Center Utrecht, Utrecht University, 3584CX Utrecht, The Netherlands	\par\medskip
	Received at Royal Society Interface 18 July 2023; Received in revised form 31 October 2023;\\ Accepted 01 November 2023; Published 22 November 2023.\\
	DOI: \url{https://doi.org/10.1098/rsif.2023.0409}}

\cortext[mycorrespondingauthor]{Corresponding author.\\
 Email address: \url{augustine.okolie@tum.de} (A.Okolie).}


\begin{abstract}
We adopt a maximum-likelihood framework to estimate parameters of a stochastic susceptible-infected-recovered (SIR) model with contact tracing on a rooted random tree. Given the number of detectees per index case, our estimator allows to determine the degree distribution of the random tree as well as the tracing probability. Since we do not discover all infectees via contact tracing, this estimation is non-trivial. To keep things simple and stable, we develop an approximation suited for realistic situations (contract tracing probability small, or the probability for the detection of index cases small). In this approximation, the only epidemiological parameter entering the estimator is $R_0$.\\
The estimator is tested in a simulation study and is furthermore applied to covid-19 contact tracing data from India. The simulation study underlines the efficiency of the method. For the empirical covid-19 data, we compare different degree distributions and perform a sensitivity analysis. We find that particularly a power-law and a negative binomial degree distribution fit the data well and that the tracing probability is rather large. The sensitivity analysis shows no strong dependency of the estimates on the reproduction number. Finally, we discuss the relevance of our findings.
\end{abstract}
	
\end{frontmatter}

\section{Introduction}
Infectious disease models have been instrumental in the study of many infectious diseases. Usually, these models are dependent on several biological parameters which can be epidemiological such as transmission, recovery, etc, or intervention parameters such as contact tracing, screening, vaccination, and others. However, most of these parameters are not or only partially known and may cause predictions from these models to lack robustness \cite{khan2020parameter} if not chosen appropriately. Missing data poses a major quantification challenge in epidemiology due to unobserved or partially observed events \cite{little2012prevention}. This makes parameter estimation essential in modelling disease spread. 
Often, the likelihood of parameters is maximised following the model predictions on sets of parameter values. In order to achieve parameter estimation, the model system property must be identifiable, i.e. estimating its parameters uniquely from the given data \citep{canto2009structural, canto2011identifiability,craciun2008identifiability}. Several estimation methods, e.g. statistically based techniques such as Approximate Bayesian computation (ABC) \cite{blum2010hiv}, Markov chain Monte Carlo (MCMC) integration \cite{o1997bayesian}, optimal control theory approach \cite{gotz2017modeling}, classical least-squares method \cite{agusto2018optimal}, and others (also see the review article~\cite{Stollenwerk2012}) have been instrumental in estimating parameters and making inferences in epidemic models. With respect to parameter estimation, contact tracing is particularly challenging, as we somehow need to estimate the fraction of contacts we miss to identify: We need to estimate something that is {\it per definitionem} unobserved.
\par\medskip

Several estimation techniques have been proposed for estimating important intervention parameters in modelling the recent covid-19 pandemic. For instance, ~\citet{manou2022estimation} obtained a best-fit model by proposing statistical methods for the underlying 
serial interval probability distribution for the covid-19 virus in Mayotte from March 2020 to January 2022. Their method 
was then used to estimate time-varying reproduction numbers 
and transmission rates observed from the collected data.\par\medskip

Only a few attempts have been made to identify parameters specific to contact tracing:
\citet{muller2007estimating} proposed a branching process approach for contact tracing in randomly mixing populations to estimate tracing probability from contact history at the onset of an epidemic, based on the theory introduced in~\citet{muller2000contact}. The derived estimator was then applied to data from contact tracing for tuberculosis and chlamydia. ~\citet{blum2010hiv} use a Bayesian framework to estimate parameters for rates of contact tracing and detection by random screening, and the method of~\citet{dyson2018targeted} is based on fitting a Yaws and Trachoma contact tracing survey data to a stochastic household model. ~\citet{Tanaka2020} took up the branching process approach to estimate the percentage of undiagnosed persons in the covid-19 pandemic with recursive full tracing. 
\par\medskip 

In this paper, we propose methods for estimating parameters in graph-based models \cite{okolie2020exact}. A stochastic SIR model on a tree-shaped contact graph is modelled such that the underlying contact structure is given by a fixed or random graph. Due to the nature of the problem, we adapt the branching process theory results for contact tracing on random trees \cite{okolie2020exact} to formulate a likelihood estimator for estimating the tracing probability and expected number of contacts. We first performed a simulation study with a Poisson degree distribution to check the performance of the maximum likelihood estimator. Thereafter, we applied the model to contact tracing data collected during the covid-19 pandemic in Karnataka, India. 
Overall, we show that our estimator based on the branching process theory for contact tracing is well suited for estimating tracing probabilities and degree distribution of the underlying contact structure in tree-based models. 
\par\medskip
The remaining part of the paper is structured as follows: Section~2 outlines the tree model and model assumptions. Section~3 presents the distribution of ages since infection, while Section~4 discusses the distribution of detected cases from one index case. We set up a likelihood estimator for estimating the tracing probability and underlying contact structure using these results and simulated data in Section~5 followed by a sensitivity analysis in Section~6. Last, we discuss our findings in Section~7.

\subsection{Related works on SARS-CoV-2 epidemic models}

Over the past few years, there has been a growing body of research on the mathematical modelling of infectious diseases, particularly the SARS-CoV-2 virus. \citet{bertacchini2020temporal} provided insights into the temporal spreading of the virus, examining key parameters that influence the rate of spread. The work of \citet{chondros2022integrated} presented an integrated simulation framework for both the prevention and mitigation of pandemics caused by airborne pathogens, providing a comprehensive approach towards understanding the dynamics of airborne diseases.
\par\medskip
Furthermore, \citet{cuevas2021lockdown} studied lockdown measures and assessed their impact on the COVID-19 outbreak in Mexico using both single- and two-age-structured epidemic models. \citet{kovacevic2020distributed}, on the other hand, employed a distributed optimal control epidemiological model for understanding the COVID-19 pandemic, emphasizing the importance of coordinated control efforts in disease mitigation.
\citet{modi2021simulation} focused on the spread of COVID-19 in India using the SEIR model, providing crucial insights into the potential dynamics of the virus in dense populations. \citet{kevrekidis2021reaction} added a spatial dimension to the modelling of COVID-19, studying the outbreak in Greece and Andalusia using reaction-diffusion models.
\par\medskip
These works, while providing valuable insights into the dynamics and control of the SARS-CoV-2 virus, largely emphasize temporal, spatial, or control aspects. Our contribution, in contrast, focuses on the intricacies of contact tracing in the context of a stochastic SIR model implemented on a rooted random tree. We emphasize the challenge posed by undetected cases and the non-trivial nature of parameter estimation in such models. By offering a unique perspective on parameter estimation and contact tracing, we aim to add depth to the existing literature and contribute to a more comprehensive understanding of disease spread dynamics as observed in the COVID-19 outbreak.

\section{Model assumptions}
For the convenience of the reader, we will first sketch the motivation and idea of the branching theory process for contact tracing on rooted random trees (tree model) in~\citet{okolie2020exact} for our estimation analysis. A contact network in most applications represents individuals as nodes and interaction links via edges. Interaction links are channels where individuals can have direct or indirect contact, e.g. family, school, work, etc. These contact networks are applicable and useful in analysing contact tracing because they hold information about individuals and their neighbours \cite{green2010large, keeling1999effects}. However, applying contact networks to infectious disease dynamics is not straightforward as it requires a detailed understanding of the underlying network structure, e.g. the degree distribution and correlation, clustering coefficients, and properties defined by the network topology.
\par\medskip
Once we have a defined contact network with pre-defined nodes and interacting links, we have a contact graph. The basic idea is to describe an epidemic by constructing a simple contact graph that is a rooted tree where only the root node is infected at the onset of the epidemic. The choice of this tree contact graph is for mathematical convenience as trees are not appropriate to describe more complex interactions for natural contact graphs. However, from a microscopic level, we can gain a better understanding of the overall mechanism and
functioning of larger and more complex graphs, as many graph models as the configuration model resemble locally a tree~\cite{kiss2017mathematics}. Then we assign independently on each edge connecting one infected and one uninfected node a probability of transmitting the disease. If we focus on edges that transmit the disease, then we have only the infection graph which is a subgraph of the contact graph. Contact tracing is also analysed on this infection graph such that upon recovery of an index case, direct neighbours of this index case are also removed with some tracing probability. From the number of detected cases by an index case via contact tracing, it is possible to estimate the degree distribution of the underlying contact network and also the tracing probability.
\par\medskip

\begin{figure}[h]
	\centering
    \includegraphics[width=.35\textwidth]{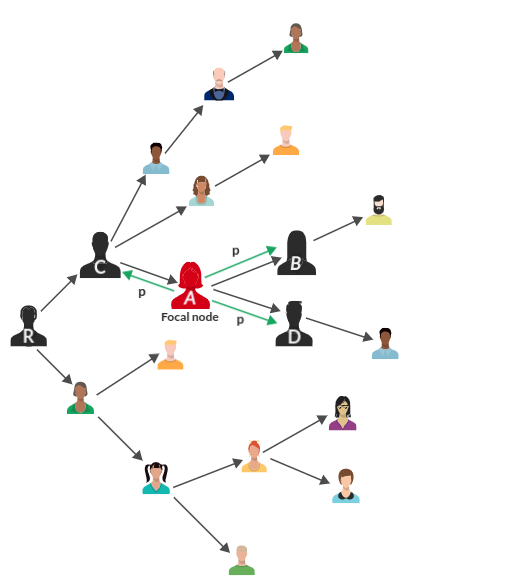}
    \caption{
    Schematic representation of the infection graph, illustrating the dynamics and interconnections in the process of disease transmission, and both forward and backward tracing. The root node is individual $R$. Focal individual $A$, infected by individual $C$, subsequently infects individuals $B$ and $D$. In a forward tracing scenario starting from individual $A$, individuals $B$ and $D$ can be traced. In contrast, in a backward tracing scenario from $A$, individual $C$ can be traced. The tracing probability is denoted in green.
    }
    \label{fig:tree}
\end{figure}

On the rooted random tree (see Fig. \ref{fig:tree}), the infection starts from the root node $R$ and spreads downwards through the directed edges. Individual $C$ which is a direct contact of the root node is infected and spreads the disease to the focal individual $A$. Furthermore, $A$ also spreads the infection to $B$ and $D$. Individuals $B$ and $D$ are the ``downstreams" (infectees) of $A$ while $C$ is an ``upstream" (infector) of $A$. We define $K$ a random variable denoting the number of downstream edges of an individual with expectation $\mathbb{E}[K]$, where we assume that the downstream degree of each node is an i.i.d.\ realisation of $K$. We note that the root node $R$ is special as it has no infector, or equivalently, it has no upstream edge. It follows that $K$ coincides with the degree only of the root, while for all other nodes have a degree of $K+1$ (infectees plus infector). At the onset of the epidemic, only the root is infected while other individuals are susceptible. We consider a susceptible-infected-recovered (SIR) model such that a recovered individual remains immune and does not get re-infected. Contact on one edge (between a susceptible and infected) will lead to infection. On a given edge, contacts happen at exponentially distributed waiting times at rate $\beta$. An infected individual recovers either unobserved at rate $\alpha$, or observed and diagnosed at rate $\sigma$. Diagnosed individuals are immediately isolated or treated and classified as recovered. With probability $p_{obs} = \sigma/(\alpha + \sigma)$, an infected individual eventually is observed. 
\par\medskip

An observed/diagnosed infected individual not only becomes isolated but also is an index case that triggers contact tracing. 
That is, every adjacent edge has an independent probability $p$ to be traced and consequently isolated if infected. In accordance with the data analysis we aim at, we focus here on one-step tracing, that is, traced individuals do not trigger further tracing events. We do, however, take into account forward and backward tracing as described in \citet{okolie2020exact}. We do note this fact, as quite often, theoretical work solely focuses on forward tracing.
\par\medskip
All in all, an infected individual can lose his or her infectivity in three possible ways; an unobserved recovery $\alpha$, observed recovery $\sigma$, and a successful tracing event. 
It turns out that the central ingredient for the analysis is the probability for an infected individual to still be infectious at age $a$. Please note that ``age'' in the present paper always refers to the age of (or time since)  infection, and never to chronological age. We define 
\begin{eqnarray}
{\kappa(a) = \mathbb{P}(\mbox{a randomly chosen infected node of generation is infectious at age of infection }a),}
\end{eqnarray}
which satisfies the following differential equation:

\begin{eqnarray*}
\frac{d}{da} \kappa(a) =  - \kappa (a)\left( {\alpha  + \sigma  + \mbox{tracing(a)} } \right),
\label{KappaSIR}
\end{eqnarray*}

where $\kappa(0) = 1$. Without contact tracing, $\hat \kappa (a) := {e^{ - (\alpha  + \sigma )a}}$. With contact tracing, this probability $\kappa(a)$ is decreased and thus

\begin{equation} \label{KappaProbInfec}
\kappa (a) = \hat \kappa (a)[1-p \times \text {tracing in the interval }[0, \mathrm{a})].
\end{equation}
In~\cite{okolie2020exact}, expressions for $\kappa(a)$ are derived. As we do not use these results in the current paper, we only indicate the overall structure and refer the interested reader to that paper for the details.

\section{Distribution of ages since infection}
In order to work out the distribution of the number of detectees per index case, the age since infection of the index case at its diagnosis is required. Thereto, we consider the case without contact tracing, $p=0$ (such that index cases are diagnosed but do not trigger contact tracing). This assumption simplifies the arguments and yields an approximation for the age distribution in case of $p>0$, which is still appropriate if $p\ll 1$ or if $p_{obs}\ll 1$. It turns out, that the resulting approximation is sufficient for practical purposes.\\ 
Since the recovery rate $\alpha$ and the screening rate $\sigma$ are constant, we have a Markovian model, and the age distribution of index cases coincides with the age distribution in the population. 

Let $i(t,a)$ denote the age since infection-structured population size of infected individuals. As derived in~\citet{okolie2020exact}, the age-structured model reads

\begin{eqnarray}
(\partial_t+\partial_a) i(t,a) &=& - (\alpha + \sigma) i(t,a)\\
i(t,0) &=& \int_0^\infty \theta(a)\, i(t,a)\, da.
\end{eqnarray}

where  
$$\theta(a)=\mathbb{E}[K]\, \beta\,e^{-\beta a}$$ 
is the age-dependent rate at which an infected average individual produces (downstream) infecteds. At this point, it is crucial that the contact graph is a tree, and the downstream degree distribution of each node/individual is an i.i.d.\ realization of the degree distribution $K$. If we count the number of nodes with a certain distance to the root, this number of nodes is exponentially increasing (for $\mathbb{E}[K]>1$), unless the tree is finite in a given realization. We can exclude the case of finite trees since these realizations imply that we have a minor outbreak, and we are not interested in these minor outbreaks. As usual, the age-structured model will tend to an exponential growing solution with a stable age structure,

$$i(t,a) = I_0\,e^{\lambda t}\,i(a)$$ 
with $i(a) = e^{-(\lambda+\alpha+\sigma) a}$ the probability to be infectious at age $a$. The exponent $\lambda$ is the unique real root of $$ 1 
= \int_0^\infty \theta(a)\, e^{-(\lambda+\alpha+\sigma)a} \, da
=  \mathbb{E}[K]\,\beta\,\int_0^\infty\,e^{-(\lambda+\alpha+\sigma+\beta)a} \, da
\qquad\Rightarrow\qquad 
\lambda = \beta (\mathbb{E}[K]-1)-\alpha - \sigma.$$ 

\begin{figure}[h!]
\centering
\includegraphics[width=0.5\textwidth]{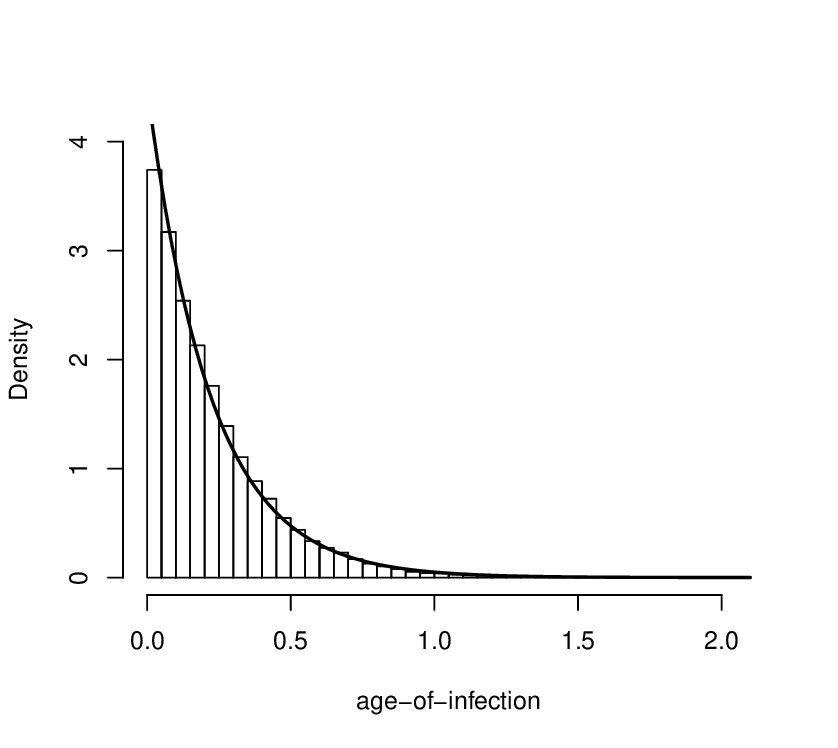}
\caption{
Theoretical age distribution $\varphi(a)$ (solid line) vs. simulated age distribution, represented by bars. Parameters: $\beta=1.5$, $\alpha=0.5$, $\sigma=0.5$, $p=0.6$, and $\mathbb{E}[K]=4$.
}

\label{fig:age-distr}
\end{figure}

The asymptotic age distribution of index cases (which are detected at rate $\sigma$) tends to
\begin{equation} \label{eq:asymAgeDistri}
\varphi (a) = \lim_{t\rightarrow\infty}\,\frac{\sigma\,i(t,a)}{\int_0^\infty \sigma\,i(t,b)\,db} = \beta\left(\mathbb{E}[K]-1\right)\,e^{-\beta\left(\mathbb{E}[K]-1 \right)a}.
\end{equation}
As a side remark, we also obtain the reproduction number
from these considerations by
$$ R_0 = \int_0^\infty \theta(a)\,  e^{-(\alpha+\sigma)a} \, da
= \frac{\mathbb{E}[K]\, \beta}{\alpha+\sigma+\beta}.$$
As shown in Fig.~\ref{fig:age-distr}, the agreement of the age distribution with simulated data is still excellent, though we have in the simulation $p=0.6$ and $p_{obs}=1/2$. Furthermore, we have a higher density of lower age groups in the population. For any randomly chosen individual given by age since infection, it is not surprising to have a younger dominating age class. Due to the exponentially fast-growing population, this asymptotic age distribution is expected.

\section{Distribution of detected cases}

In this section, we derive the distribution of the number of detected cases per index case. That is the fraction of contacts of an index case who are detected via a tracing event triggered by the index case. We start with forward tracing. We then combine this result with backward tracing to yield full tracing.\\
Note that a central ingredient is the age distribution derived in the last section. We did not include contact tracing there. That is, all results in the present section are only a valid approximation if contact tracing does not crucially affect this age distribution. This is the case if either $p$ or $p_{obs}$  is small. All results are only valid under this assumption. However, the simulation study discussed below shows that this assumption is not too restrictive for practical purposes.

\begin{prop}
Let $\hat{p}(a)$ be the probability that an infected downstream node is successfully traced given that the focal individual becomes an index case at age since infection $a$. 

\begin{eqnarray}
\hat{p}(a) = p\frac{\beta}{\alpha + \sigma - \beta} \big(e^{-\beta a} - e^{-(\alpha + \sigma) a} \big)
\end{eqnarray}

\end{prop}

\begin{proof}

Note that an individual is only able to become an index case at the transition from I to R, that is, our focal individual is infectious in $[0,a)$.  We consider one downstream individual. Let $s_1(c)$ represent the probability for this downstream individual to still be susceptible at age $c\in[0,a]$, $s_2(c)$ the probability to be infected, and $s_3(c)$ the probability  
for the downstream node to be removed. 
We have the following ODEs

\begin{eqnarray*}
\dot{s_1} &=& - \beta\,s_1 \qquad \qquad \qquad s_1(0) = 1 \\
\dot{s_2} &=& \beta s_1 - (\alpha + \sigma)\,s_2 \quad \,\,\,\,\, s_2(0) = 0 \\
\dot{s_3} &=& (\alpha + \sigma)\,s_2 \qquad \qquad \,\, s_3(0) = 0.
\end{eqnarray*}

\begin{figure}[h!]
\centering
\includegraphics[width=0.5\textwidth]{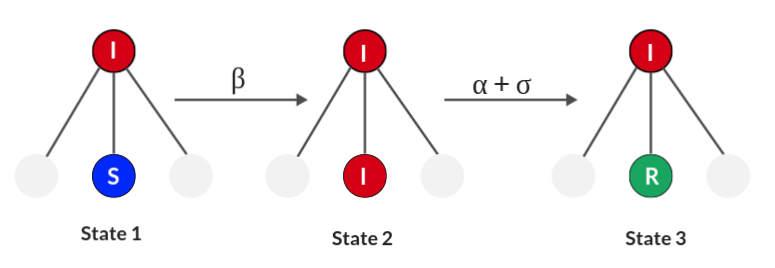}
\caption{
Transition states for a single edge in the infection process. These states represent the probabilities $s_1(c)$, $s_2(c)$, and $s_3(c)$ of an individual downstream from the index case to remain susceptible, infected, and be removed at age $c \in [0,a]$.}
\label{fig:p_state}
\end{figure}

The probability for the downstream node to be infectious at the time the infector has age since infection $a$ given by $s_2(a)$,  
$$
s_2 = \frac{\beta}{\alpha + \sigma - \beta} \big(e^{-\beta a} - e^{-(\alpha + \sigma) a} \big),
$$
and $\hat p(a) = p s_2(a)$ establishes the result. 
\end{proof}
With this proposition and the age distribution $\varphi(a)$, we are able to find the distribution of the number of detected downstream individuals. For simplicity, we first consider a fixed degree distribution $K=k$ for some $k\in\N$, and then address the case of a random tree, where $K$ is a random variable.

\subsection{Fixed degree}
In the present section,  assume that the downstream degree of a node in the tree always is a deterministic number $K=k\in\N$. Particularly, $\mathbb{E}[K]=K$.
\begin{prop}
Let $T$ be the random variable for the total number of successfully traced individuals by one index case and forward tracing only. The asymptotic probability distribution of $T$ under forward tracing reads
\begin{eqnarray}\label{T_result}
P(T = i) = \int_{0}^{\infty}{k \choose i}\,\hat{p}(a)^{i}\,\left(1- \hat{p}(a) \right)^{k-i}\beta\left( k -1 \right)\,e^{-\beta \big( k-1 \big)a}\,da.
\end{eqnarray}
\end{prop}

\begin{proof}
As we assume that a tracing event acts independently on different edges, the random variable $T$, conditioned on the age of the index case at diagnosis $a$,  follows a Binomial distribution with parameters $k$ and an age-dependent tracing probability on one edge $\hat{p}(a)$, $T \sim Binom(k, \hat{p}(a))$. Thus, the probability of $i$ downstream detectees given age $a$ and $k$ total downstream nodes is given as 

\begin{eqnarray}
\label{eq:fixed}
P(T=i\,|\,a) = {k \choose i}\,\hat{p}(a)^{i}\left( 1 - \hat{p}(a)\right)^{k-i}.
\end{eqnarray}

Last we remove the condition $a$ by integrating over all possible age of index cases $\varphi(a)$ (eqn. \ref{eq:asymAgeDistri}),

\begin{eqnarray}
P(T=i) &=& \int_{0}^{\infty}P(T=i\,|\,a)\,\varphi(a)\,da 
 =  \int_{0}^{\infty}{k \choose i}\,\hat{p}(a)^{i}\,\left(1- \hat{p}(a) \right)^{k-i}\beta\left(k -1 \right)\,e^{-\beta \big( k-1 \big)a}\,da
\label{eqn:detecte_given_a}
\end{eqnarray}

\end{proof}

Now we turn to full tracing. Thereto, we introduce the random characteristic $I_a$, which assumes the value $1$ if the upstream individual of the index case (its infector) still 
is infected when the index case is identified (where the index case has age since infection $a$), and $0$ else. Note that 
 $I_a$ is a Bernoulli random variable with 
 $$ P(I_a=1) = e^{-(\alpha + \sigma)a} + \mathcal{O}(p).$$
As before, in what follows we use the approximation
 $$ P(I_a=1) = e^{-(\alpha + \sigma)a}$$
 and drop the ${\cal O}(p)$ correction terms.
 
\begin{prop}
Let $T_{tot}$ be the random variable for the total number of successfully traced individuals by one index case, under full tracing (forward and backward tracing). With $\varphi(a)$ and $I_a$ as introduced above, the probability distribution of $T_{tot}$ reads
\begin{eqnarray}\label{eq:final_eqn}
P(T_{tot} = i) = \int_{0}^{\infty} \bigg[pP(I_{a} = 1)\,P(T=i-1\,|\,a) + \big(1 - pP(I_{a} = 1)\big)\,P(T=i\,|\,a) \bigg]  \varphi(a)\,da.
\end{eqnarray}
\end{prop}

\begin{proof}
If the infector already is recovered ($I_a=0$), then 
(conditioning on the age/time of infection of the index case $a$)
$$ P(T_{tot}=i|a,\,\, I_a=0) = P(T=i|a).$$
If the infector is still infectious, also the infector might be traced, such that one of the $k$ detectees might be the upstream individual (probability $p$), or not (probability $1-p$), 
$$ P(T_{tot}=i|a,\,\, I_a=1) = p\,P(T=i-1|a)+ (1-p)\,P(T=i|a).$$
Taking these two cases together, we have 
\begin{eqnarray*}
P(T_{tot}=i|a) 
&=& 
P(I_a=0)\, P(T=i|a)
+
P(I_a=1) \bigg(\, p\,P(T=i-1|a)+ (1-p)\,P(T=i|a)\, \bigg)\\
&=& pP(I_{a} = 1)\,P(T=i-1\,|\,a) + \big(1 - pP(I_{a} = 1)\big)\,P(T=i\,|\,a).
\end{eqnarray*}
Integrating by $\varphi(a)\, da$ removed the condition on $a$ and yields the result.
\end{proof}

\subsection{Random degree}

\begin{figure}[h!]
\centering
\includegraphics[width=0.49\textwidth]{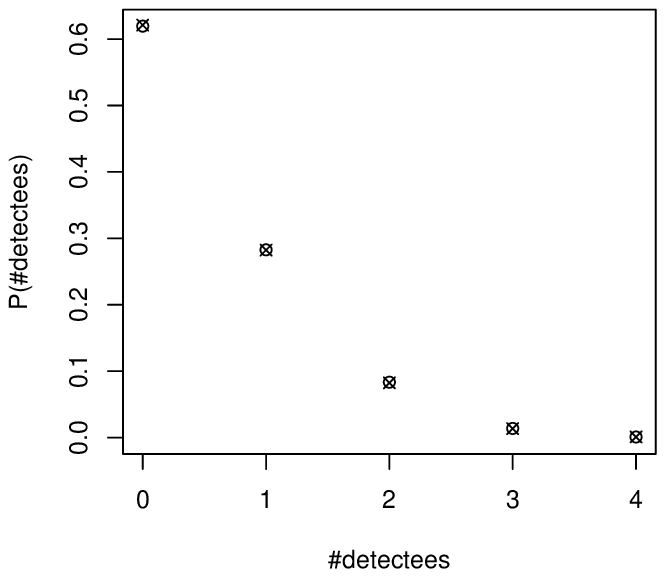}
\includegraphics[width=0.49\textwidth]{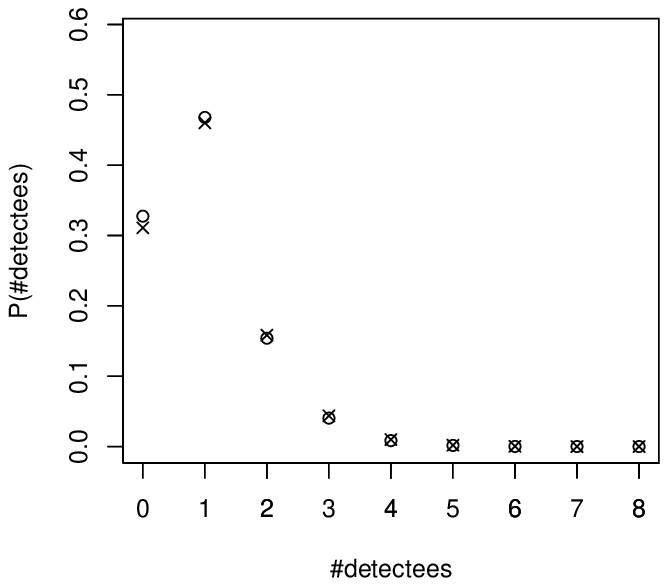}
\caption{
Distribution of detected cases per index case: theoretical predictions vs. simulated results. The left panel represents forward tracing for fixed degrees, while the right panel depicts full tracing for a Poisson graph. In both panels, crosses denote theoretical probabilities $P(T=i)$ or $P(T_{tot}=i)$, whereas circles represent results from 100000 simulations. Additional parameters used during the simulation: $\beta=1.5$, $\alpha=0.5$, $\sigma=0.5$, $p=0.6$, and $\mathbb{E}[K]=4$.
}
\label{fig:detectee-distr}
\end{figure}

So far, the model is formulated for a fixed degree. In most applications, we do not always know individual contacts $k$ due to randomness in contact structure. We now assume an arbitrary degree distribution such that the distribution of the contacts of a random node is defined by some probability distribution  $P(K=k)$. The model for fixed case in eqn.~\ref{eq:fixed} is adapted, we only have to take the expectation by summing over all possible numbers of contacts $k$ multiplied by the corresponding probabilities. Thus,  

\begin{eqnarray}
\label{eq:final_eqn_random}
P(T_{tot}=i) = 
\sum_{k=i}^{\infty}
\int_{0}^{\infty} \bigg[pP(I_{a} = 1)\,P(T=i-1\,|\,a,\,\,K=k) + \big(1 - pP(I_{a} = 1)\big)\,P(T=i\,|\,a,\,\,K=k) \bigg]  \varphi(a)\,da\,  P(K=k).
\end{eqnarray}
As illustrated in Fig. \ref{fig:detectee-distr}, we have the distribution of the number of detected secondary cases via contact tracing. For the parameter chosen in our study (see Fig. \ref{fig:k-p-distr}), we find a good agreement between our theory results and simulation. We again emphasize that the age structure entering our estimation only is an approximation, as contact tracing is neglected. nevertheless, the results are more than acceptable, even for $p=0.6$ and $p_{obs}=0.5$. 

We note that our estimator is independent of time $t$. The only ``time'' that appears is the age since infection $a$. The probability $P(T_{tot}=i)$ consists of integrals as 
$ \int_0^\infty g(a)\,\varphi(a)\, da$. 
Here, we are allowed to choose the time unit, resp.\ to define 
$a = \zeta b$ for $\zeta>0$, 
$$ \int_0^\infty g(a)\,\varphi(a)\, da
=
 \int_0^\infty g(b\, \zeta)\,\zeta \varphi(b\,\zeta)\, db.$$
If we choose one time unit to be $1/(\alpha+\sigma)$, which is, 
$\zeta=1/(\alpha+\sigma)$, then in in $P(T_{tot}=i)$ the epidemiological parameters $\beta$, $\alpha$ and $\sigma$ always can replaced by an expression of $R_0$ and $E[K]$. That is, the epidemiological parameters enter the estimator solely via $R_0$. We only check that fact for one of the terms, as the argument is similar for the other terms.
\begin{eqnarray*}
&&\int_0^\infty P(I_a=1)\,\, P(T=i\,|\, a, K=k)\,\varphi(a)\, da
\\
&=&{k\choose i} 
\int_0^\infty e^{-(\alpha + \sigma)a}\,\, \hat{p}(a)^{i}\,\left(1- \hat{p}(a) \right)^{k-i}\,
\,\beta\left(\mathbb{E}[K]-1\right)\,e^{-\beta\left(\mathbb{E}[K]-1 \right)a}\, da
\\
&=&{k\choose i} 
\int_0^\infty e^{-b}\,\, \hat{p}(b/(\alpha+\sigma))^{i}\,\left(1- \hat{p}(a/(\alpha+\sigma)) \right)^{k-i}\,\frac{\beta\left(\mathbb{E}[K]-1\right)}{\alpha+\sigma}\,e^{-\frac{\beta\left(\mathbb{E}[K]-1 \right)}{\alpha+\sigma}\,b}\, db.
\end{eqnarray*}
With
$ R_0 = \frac{\beta\,\mathbb{E}[K]}{\alpha+\sigma+\beta}$
we have (note that always $\mathbb{E}[K] > R_0$, as we have -- in average -- only $\mathbb{E}[K]$ downstream individuals who can get infected)
$$ \frac{\beta}{\mu+\sigma} = \frac{R_0}{\mathbb{E}[K]-R_0}$$
and hence
$$ \frac{\beta\,(\mathbb{E}[K]-1)}{\alpha+\sigma} = 
\frac{(\mathbb{E}[K]-1)\, R_0}{\mathbb{E}[K]-R_0}, \quad 
\hat p(b/(\alpha+\sigma)) = 
\frac{\frac{\beta}{\alpha+\sigma}}{\frac{\beta}{\alpha+\sigma}-1}\, 
\bigg(e^{-b}-e^{-\frac{\beta}{\alpha+\sigma}\, b}\bigg)
= 
\frac{R_0}{2\,R_0-\mathbb{E}[K]}\, \bigg(e^{-b}-e^{-\frac{R_0}{\mathbb{E}[K]-R_0}\, b}\bigg).
$$
That is, all expressions only depend on $R_0$, $K$, and $p$. 
\begin{cor}
$P(T_{tot}=i)$ and $P(T=i)$  only depend on the epidemiological parameters via $R_0$ and depends furthermore on the degree distribution given by $K$ and on the tracing probability $p$.
\end{cor}
We can use the formulas from above, where we pragmatically set $\alpha+\sigma$ to $1$, and -- given $R_0$ and $\mathbb{E}[K]$ -- define $\beta = R_0/(\mathbb{E}[K]-R_0)$.

\section{Estimation by maximum likelihood method}
In this section, we will set up the likelihood estimator for our model. We assume that we have $n$ observations of index cases, and where $i_{\ell} \in \mathbb{N}_{0}, \ell=1, ., n$ denote the total number of detectees per index case (one step tracing only). 

\subsection{Likelihood estimator}
We are able to set an estimator via $P(T_{tot}=i) = P(T_{tot}=i \mid \vec{\mu})$ for these data points where $\vec{\mu}$ are the parameters of the model we wish to estimate (tracing probability and parameters of the random variable $K$, e.g.\ $\mu=(p, \mathbb{E}[K])$ in case of a Poisson distribution for $K$). The likelihood for the data reads
\begin{eqnarray*}
\mathcal{L}\left(\vec{\mu} \mid i_{\ell}, \ell=1, . ., n \right)=\prod_{\ell=1}^{n} \int_{0}^{\infty} \bigg[pP(I_{a} = 1)\,P(T=i_l-1\,|\,a) + \big(1 - pP(I_{a} = 1)\big)\,P(T=i_l\,|\,a) \bigg]  \varphi(a)\,da
\label{Eqn:likelihood}
\end{eqnarray*}

and the log-likelihood is given by 
\begin{eqnarray}
\small
\mathcal{L} \mathcal{L} \left(\vec{\mu} \mid i_{\ell}, \ell=1, ,.., n\right)= \sum_{\ell=1}^{n} \ln \bigg( \int_{0}^{\infty} \bigg[pP(I_{a} = 1)\,P(T=i_l-1\,|\,a) + \big(1 - pP(I_{a} = 1)\big)\,P(T=i_l\,|\,a) \bigg]  \varphi(a)\,da \bigg).
\label{Eqn:loglikelihood}
\end{eqnarray}

\subsection{Simulated data}
\begin{figure}[h!]
\centering
\includegraphics[width=0.42\textwidth]{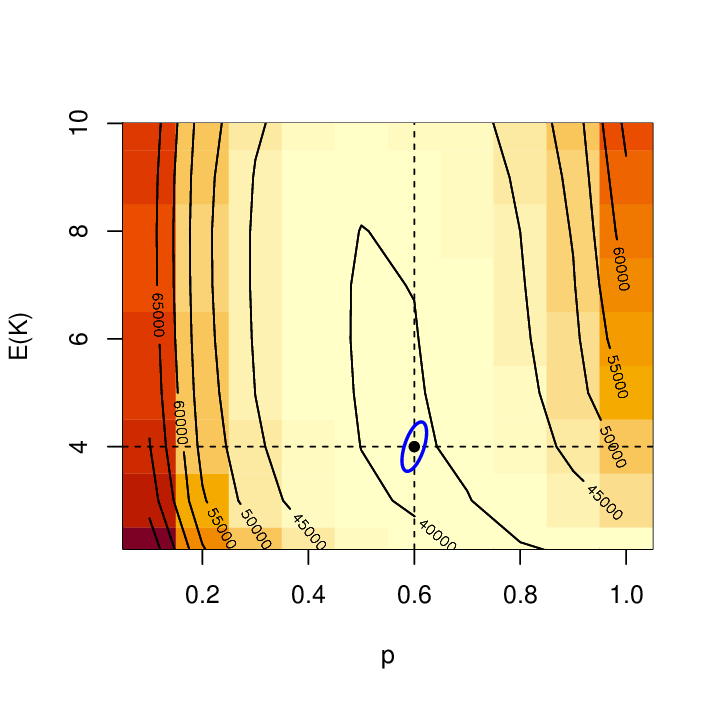}
\caption{Parameters estimated from the simulation data derived from 100000 iterations. The tracing probability, denoted as $p$, was determined to be 0.6 and the expected number of edges, represented by $\mathbb{E}[K]$, was 4. The blue region represents the 95\% confidence region for the estimated parameters. Additional parameters used during the simulation: $\beta=1.5,\,\alpha=0.5,\,\sigma=0.5$.}
\label{fig:k-p-distr}
\end{figure}

An agent-based stochastic model is used to simulate the data. To maximise the likelihood in eqn. \ref{Eqn:loglikelihood}, we plug into the likelihood function all independent observed data points and determine the arg max. As shown in Fig. \ref{fig:k-p-distr}, in our estimation we find back the true values we used in the simulation, namely the tracing probability $p=0.6$ and the expected number of edges $\mathbb{E}[K]=4$. 
The blue circle region contains a global maximum for the true parameter value of the estimated Poisson degree distribution. Other parameters $\beta,\,\alpha,\,\sigma$ are known and fixed. For the simulated data, we find a satisfying result based on our theory assumption.

\subsubsection{Stability of the estimator against non-tree typologies in the contact graph}

Real-world networks are, of course, no trees. We investigate the stability of our estimators against the violation of this prerequisite. Our estimator is particularly based on two assumptions: We know the distribution of the index cases' time since infection, and all downstream nodes are susceptible when a node becomes infected.\\ 
The notion of downstream nodes as introduced above depends on the tree topology. We generalize this notion in calling all neighbours of a focal infected node ``downstream nodes'' apart from the infector of this focal node. Circles and clusters in a contact graph might lead to infections of downstream nodes from outside (which we call ``outside infections''). \\
The deviation from a tree can be measured on the topological level or on an immunological level. Speaking about topology, particularly the appearance of triangles is well known to affect theory which is based on a tree topology: The message-passing method~\cite{newman2002spread,karrer2010message} is an exact version of the pair approximation on trees. Pair approximation, in turn, requires correction terms if considered on more general networks~\cite{keeling1999effects}. We thus expect that triangles might challenge our estimator. We use the configuration model~\cite{kiss2017mathematics} in an adapted version which allows to control the fraction of nodes in triangles (see supplementary information). \\
The second level where the deviation from trees becomes visible is the epidemiological process: The fraction of outside infections is an alternative characterization for non-tree graphs. It turns out (see supplementary information) that this epidemiological characterization better predicts the performance of the estimator than the density of triangles: Outside infections slow down the spread of an epidemic, and in that, the time-since-infection structure of index cases is shifted to longer infectious periods. In that, index cases have more time to infect their downstream nodes. Moreover, outside infections produce potentially even more infected downstream nodes as we expected from tree-based models. Both effects point in the same direction, such that $p$ tends to be overestimated to explain the additional detectees, which are induced by the epidemiological consequences of the graph topology.\par\medskip 

\begin{figure}[htb]
\begin{center}
(a)\includegraphics[width=5cm]{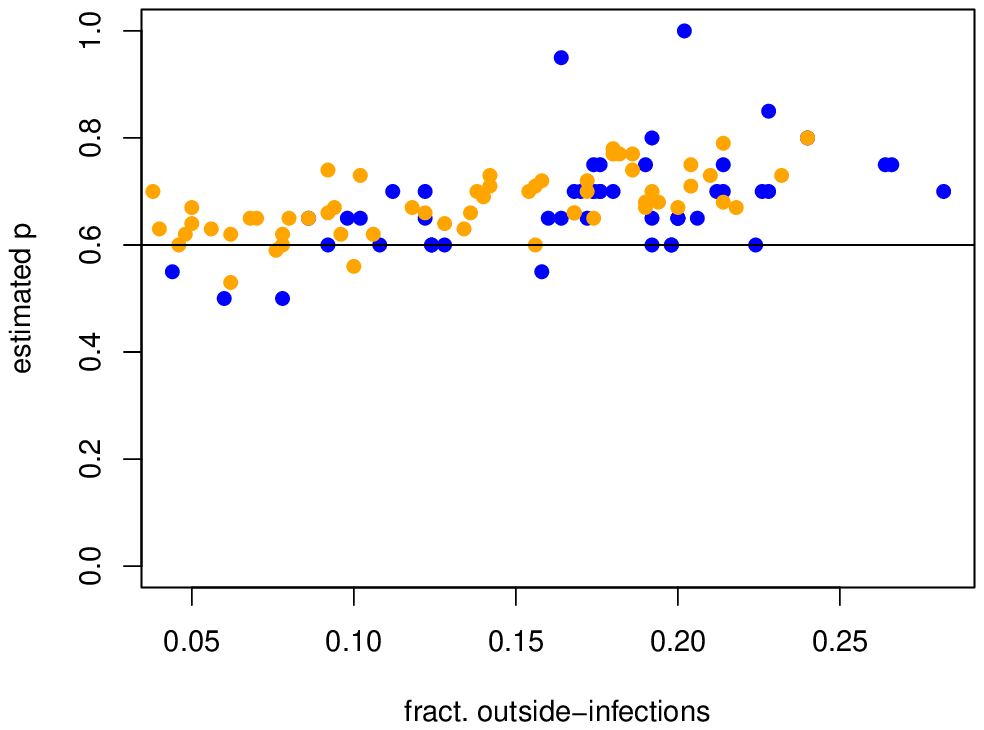}
(b)\includegraphics[width=5cm]{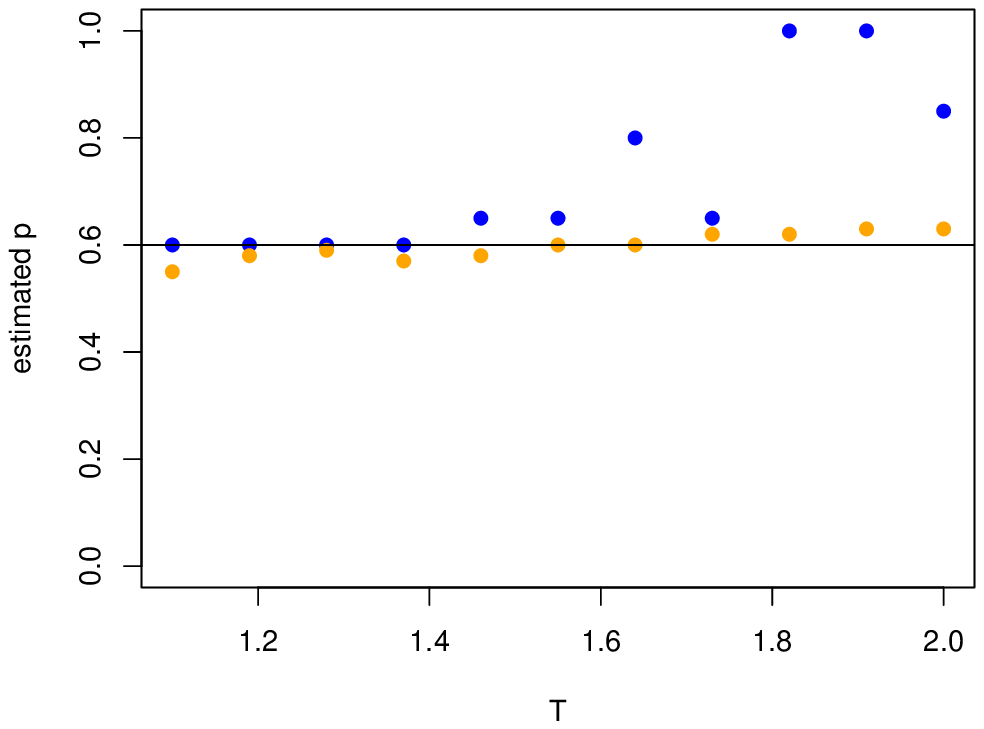}
(c)\includegraphics[width=5cm]{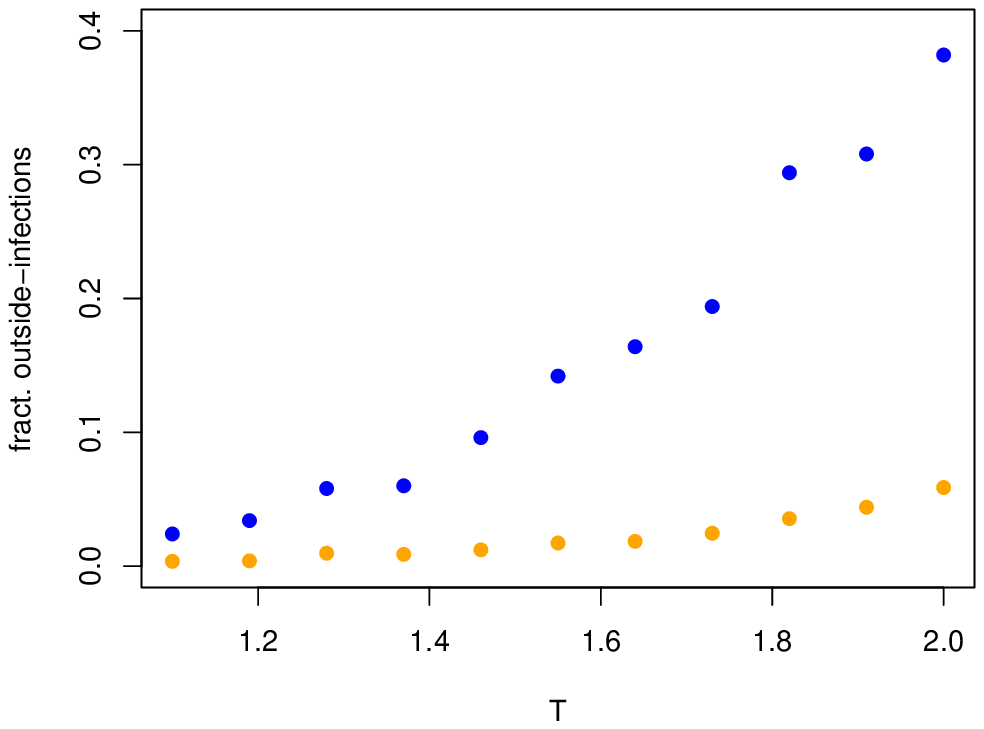}
\end{center}
\caption{Results of the simulation study for the performance of the estimator for the generalized configuration model (orange: fixed excess degree, blue. Poisson excess degree). (a) Estimated $p$ over the fraction of index cases with outside infections, induced by forcing the configuration graph to have triangles. For the fixed degree, index cases from the time interval $[0,2]$ are used, for the Poisson degree distribution, we use the time interval $[0,1.2]$. (b), (c) We consider the index cases identified in the time interval $[0, T]$, where $T$ is on the x-axis. The standard configuration model is used to produce the graph, without forcing for additional triangles. (b) Estimation of $p$. (c) Fraction of index cases with outside-infections. Parameters: $\beta=1.5$, $\alpha=\sigma=0.5$, $p=0.6$, $\mathbb{E}(K)=4$, only forward tracing.}\label{simulStudy}
\end{figure}

If we inspect the simulation study (Figure~\ref{simulStudy}, panel(a)), we find that indeed the estimates of $p$ increase over the fraction of index cases which possess outside infections, which are introduced by triangles. However, up to a fraction of 10-20\%, this effect is not too severe. \\
The next, interesting question is the influence of the excess degree distribution instead of the influence of triangles on the fraction of outside-infections. If we compare a configuration model with a fixed excess degree and a Poisson excess degree over a first time interval $[0, T]$, we find in the Poisson graph a much faster increase of these outside-infections in $T$. The difference between fixed and Poisson excess degree is much more important than the number of triangles (see supplementary information). Indeed, for the Poisson excess degree, the estimator becomes biased even during the late exponential growing phase, while in the fixed degree model, the estimations are acceptable basically during the complete exponential growing phase (Figure~\ref{simulStudy}, panel(b)); the reason is the striking difference in the number of index cases with outside-infections (Figure~\ref{simulStudy}, panel(c)). The explanation is well known: Configuration models locally look like trees. However, as it is known from the celebrated friendship paradox~\cite{Feld1991,allen2008mathematical, Christakis2010}, the preferential mixing of the configuration model lets nodes connect to nodes with a high degree. Therefore, we find in the case of the Poisson excess degree a relatively clustered, small subgroup of nodes, where the infection will take place first. We have a kind of well connected core group, which is distinctly smaller than the population size. Thus, outside infections are likely to take place after the infection moves into this core group. In the fixed degree, we of course cannot have such a core group, which explains the stability of the estimations in the fixed degree model.\par\medskip 

In the long run, however, the assumptions of an exponentially growing prevalence will not be given any more, independent of the excess degree distribution. The prerequisites of our estimator are not valid any more. In section 4 of supplementary information, we indicate a possibility of extending the basic ideas developed for trees, in order to also cover the long term behaviour of an epidemic. However, this question is not the focus of this present work.\par\medskip 

A central question now is the applicability of the theory to real world data on the background of these simulation results. This might depend on the transmission mechanisms. Sexually transmitted diseases (std's) are known to depend on core groups, and the contact network is less dynamic as in respiratory diseases. That is, we should be careful in applying our estimator to std's. The infectious contacts of respiratory infections are known to exhibit over-dispersion~\cite{lloyd2005superspreading}. However, as many of the contacts are rather casual, the infection network is not static and we won't find a fixed core group. In that, we expect that the tree-based estimators will work fine for respiratory infections.

\subsection{Covid-19 data}
In the previous section, we used the maximum likelihood estimator with contact tracing on simulated data to estimate the tracing probability and expected number of edges. In the simulated dataset, we analysed, we have information about the total number of contacts of index cases and also the number of infected contacts identified via contact tracing. In this section, we would like to see how this method works with empirical data. We wish to estimate the tracing probability $p$ and expected number of edges $\mathbb{E}[K]$ from a data set collected during the Covid-19 pandemic. We obtained a published data set on contact tracing conducted in 2020 from a remarkable extensive and nice study in Karnataka, India where 956 cases with confirmed forward contact tracing were reported between the 9th of March 2020 and the 20th of May, 2020~\citep{gupta2022contact}. A comprehensive description of the dataset including the data source, data handling and ethics approval can be found in~\citet{gupta2022contact}. A summary of the number of detectees per index case from the reported dataset is shown in Table \ref{tab:covid_data}. 

\par\medskip
We only look at one-step tracing at the moment because only detected secondary cases of primary (index) cases are accounted for in the likelihood estimator. Generally, in most epidemic modelling studies, secondary cases are defined as close/direct contacts (e.g. household, family, etc.) of index cases \cite{bell1994multistate}. In the dataset, direct and indirect detected contacts of index cases were reported. For our study, we would focus on only the direct contacts, thus we are able to analyse this scenario with the estimator for only one-step forward tracing.  For the reported reproduction number, we chose $R_0=3$ in accordance with \citet{gupta2022contact}. An extensive meta-analysis of covid-19 data from China encompassing 29 studies revealed an approximate $R_0$ value of 3.32 (95\% CI: 2.81-3.82) \cite{locatelli2021estimating}. In order to investigate the influence of $R_0$ on our estimates, we additionally carried out a sensitivity analysis. 
\par\medskip

\begin{table}[h!]
\small
\centering
\begin{tabular}{|cccccccccccccccccccc|}
\hline 
Number of detectees & 0 & 1 & 2 & 3 & 4 & 5 & 6 & 7 & 8 & 10 & 11 & 12 & 13 & 15 & 16 & 19 & 22 & 28 & 29 \\ 
\hline 
Frequency & 766 & 87 & 34 & 19 & 16 & 12 &  3 &  4 &  3 & 2 & 1 & 1 & 1  & 1 & 1 & 2 & 1 & 1 & 1 \\
\hline 
\end{tabular} 
\caption{Total number and frequencies of the total number of detected cases.}
\label{tab:covid_data}
\end{table}

We use standard random graph models (see Table \ref{tab:random_graphs}) to get inspiration on which degree distribution might be appropriate for describing the data \cite{kiss2017mathematics}.

\begin{table}
\centering
\caption{Examples of standard random graph models.}
\label{tab:random_graphs}
\begin{tabular}{ll}
\hline
\textbf{Network Model} & \textbf{Degree Distribution for Large Population Size $N$}\\
\hline
Full Graph/Random Mixing\footnotemark & $K=N-1\rightarrow\infty$, $\beta\rightarrow 0$, $R_0$ constant\\
Erd\"os-R\'enyi & $K\sim$ Poisson\\
Configuration Model & Choice: $K\sim$ Geometric\\
Scale-free Network & Power-Law, $P(K=k)=c\,k^{-\gamma}$, $\gamma>1$\\
Standard Degree Distribution & Negative Binomial\\
\hline
\end{tabular}
\end{table}

\footnotetext{See Appendix \ref{randMixAppend} and \ref{Optim} for further detail on full graph/random mixing and the optimization process respectively.}

We performed (a) a maximum likelihood estimation and (b) for model comparison, we 
used the AIC and a chi-square (goodness-of-fit) test. The summary of the point estimates is shown in Table \ref{tab:updated_table}.\par\medskip  
\textbf{Optimization}: We inspected the gradient of the result and the eigenvalues of the Hessian 
to ensure that we have (at least) an approximate local maximum. The estimator converged satisfyingly for all 
models except for the Poisson degree model: In that case, the parameter of the distribution (the expectation) always increased. Seemingly, the optimum is either very large or even infinite.\par\medskip 

\begin{figure}[tbh]
    \centering
    \includegraphics[width=\textwidth]{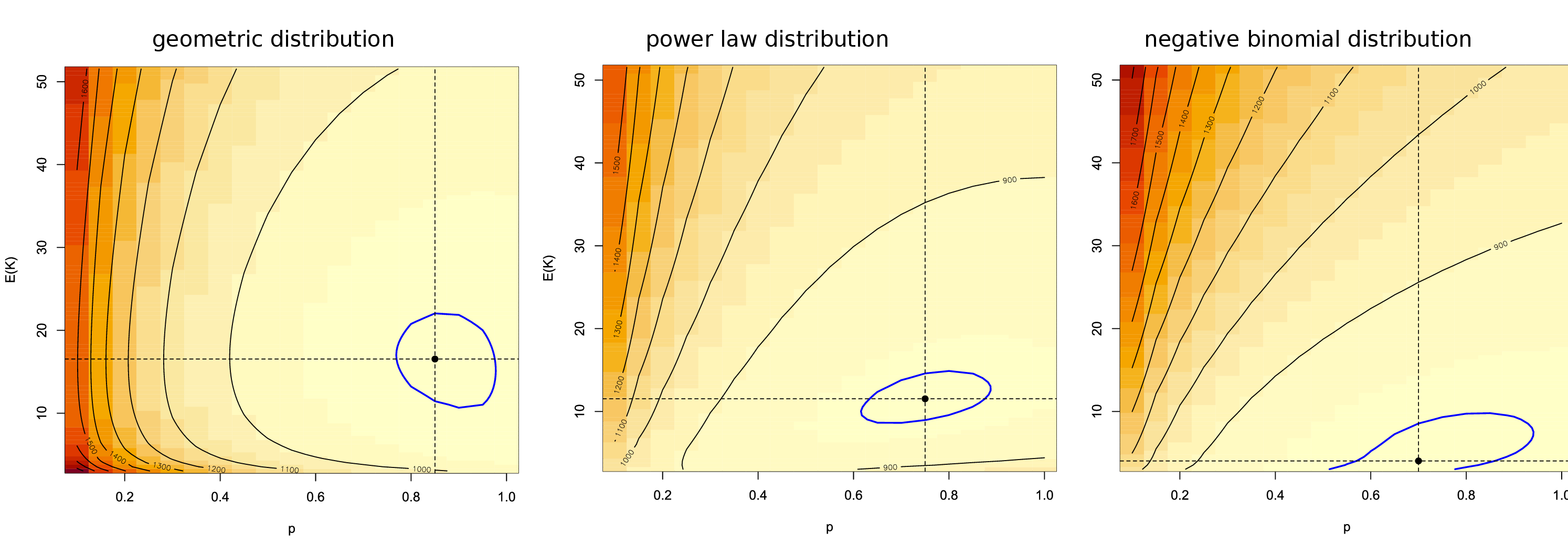}
    \caption{A contour plot showing the point estimates $p$ and $E(K)$ of the likelihood estimator in different network models. Distributions from left to right: Geometric, scale-free, and negative binomial. Parameters; $R_{0}=3$, for the negative binomial, $r=0.17$.}
    \label{fig:contour_plot}
\end{figure}

\textbf{Confidence intervals}: The approximate 95\% confidence intervals are based on the quadratic approximation of the log-likelihood at its maximum, respectively the approximation of the Fisher information matrix by the the inverse of the negative Hessian. We determined the confidence intervals for geometric, scale-free, and negative binomial distributions, as the other models turn out to be inappropriate for the data. In the case of the negative binomial distribution, however, the log-likelihood attains its maximum close to the theoretical lower boundary of $\mathbb{E}[K]$ which is $R_0$. As this is numerically a delicate situation, we approximated the confidence interval for $\mathbb{E}[K]$ in that case by the maximum value of $\mathbb{E}[K]$ such that the log-likelihood, given the shape parameter, is larger than its maximum minus $2$ (also see the blue curve in Fig.~\ref{fig:contour_plot}).\par\medskip 
\textbf{Chi square (goodness-of-fit) test}: We bin the index cases with 5-7 detectees, and all index cases with 
more than 7 detectees to ensure that at least 10 observations are in one class.\par\medskip 

\begin{table}[]
\centering
\begin{tabular}{lllllll}
\hline
Distribution & $p$ (95\% CI) & $\mathbb{E}[K]$ (95\% CI) & add. information & AIC & p-val ($\chi^2$)\\
\hline
random mixing & 0.98 & - &   & 2443 & $<10^{-20}$ & - \\
Poisson & 0.98 & 52 & (not converged) &2464 &  $<10^{-20}$\\
Geometric & 0.87 (0.76, 0.97) & 16.6 (11.3, 22.0) & & 1858 &  $<10^{-20}$\\
Power-law & 0.74 (0.61, 0.88) & 11.4 (8.3, 14.5) & $\gamma=1.48$ & 1687 & 0.003\\
negative Binomomial & 0.72 (0.59, 0.86) & 4.5 (3, 9.5) & $r=0.16$ (0.12, 0.20) & 1675 & 0.14\\
\hline
\end{tabular}
\caption{Parameter estimates and model comparison for five probability distributions fitted to a dataset, including the probability ($p$), the expected value ($\mathbb{E}[K]$), additional information (if available), Akaike Information Criterion (AIC), and p-value from a chi-squared test. Note that the estimator in the case of the Poisson distribution did not converge, and we simply fixed a large expected value. The interval for $r$ in the negative Binomial distribution is 95\% CI.}
\label{tab:updated_table}
\end{table}

We find that random mixing graphs and the Erd\"os-Renyi graph induce degree distributions with a rather lightweight tail. Therefore, these distributions do not fit the data appropriately. In the Poisson distribution, which is the degree distribution of an Erd\"osch-Renyi  graph, the optimization routine even does not 
obtain a local maximum: It seems as in this case, the optimum is only assumed for an expectation 
$\mathbb{E}[K]$ that is unreasonably large or even tends to infinity.\\
The geometric distribution has a tail that is heavy enough to at least allow for a reasonable fit, 
but AIC as well as the chi-square (goodness of fit) test indicate that this model is rejected. The power-law (scale-free graph) is the first model that is at least weakly in line with the data: The tail of a power-law distribution may become heavy, and in this, there is a possibility to handle superspreading events appropriately. The AIC is worse but not too far from the winning model, and the $p$-value for the  chi-square (goodness-of-fit) test is at least only in the range of $10^{-3}$, and not less than $10^{-20}$, as in the previous models.\\
The best model clearly is the negative binomial distribution, which is known to be an appropriate model for the number of contacts relevant to the transmission of respiratory infections~\cite{mossong2008social}. The expected number of infectious contacts is small enough to be in the range of $R_0$, while -- as expected -- the 
over-dispersion is distinct. This model has the best AIC among all models, and the goodness-of-fit test does not reject this model. \\
It is interesting that the point estimate for $p$ decreases if we choose models that have more mass in the tail. The reason is that the probability for $k$ detectees scales with $p^k$, such that $k$ small(er) 
needs to be balanced with more probability mass in the tail. A similar reason leads to smaller values for $\mathbb{E}[K]$ if the model distribution has more probability mass in the tail. However, this point estimate is rather similar for the power-law and the negative binomial and also does not heavily depend on the choice of $R_0$ (see ``sensitivity analysis''). In that, the range of $p$-estimate seems to be trustworthy. Moreover, the information in the data is sufficiently strong to point to a specific degree distribution (negative binomial),  which was not clear from the beginning. As the data are rather simple, the information content could also have been too little to allow for distinct conclusions. That the negative binomial distribution, which is well known to be appropriate in this situation, is selected, is another sign that the estimates are trustworthy.

\begin{figure}[tbh]
    \centering
    \includegraphics[width=\textwidth]{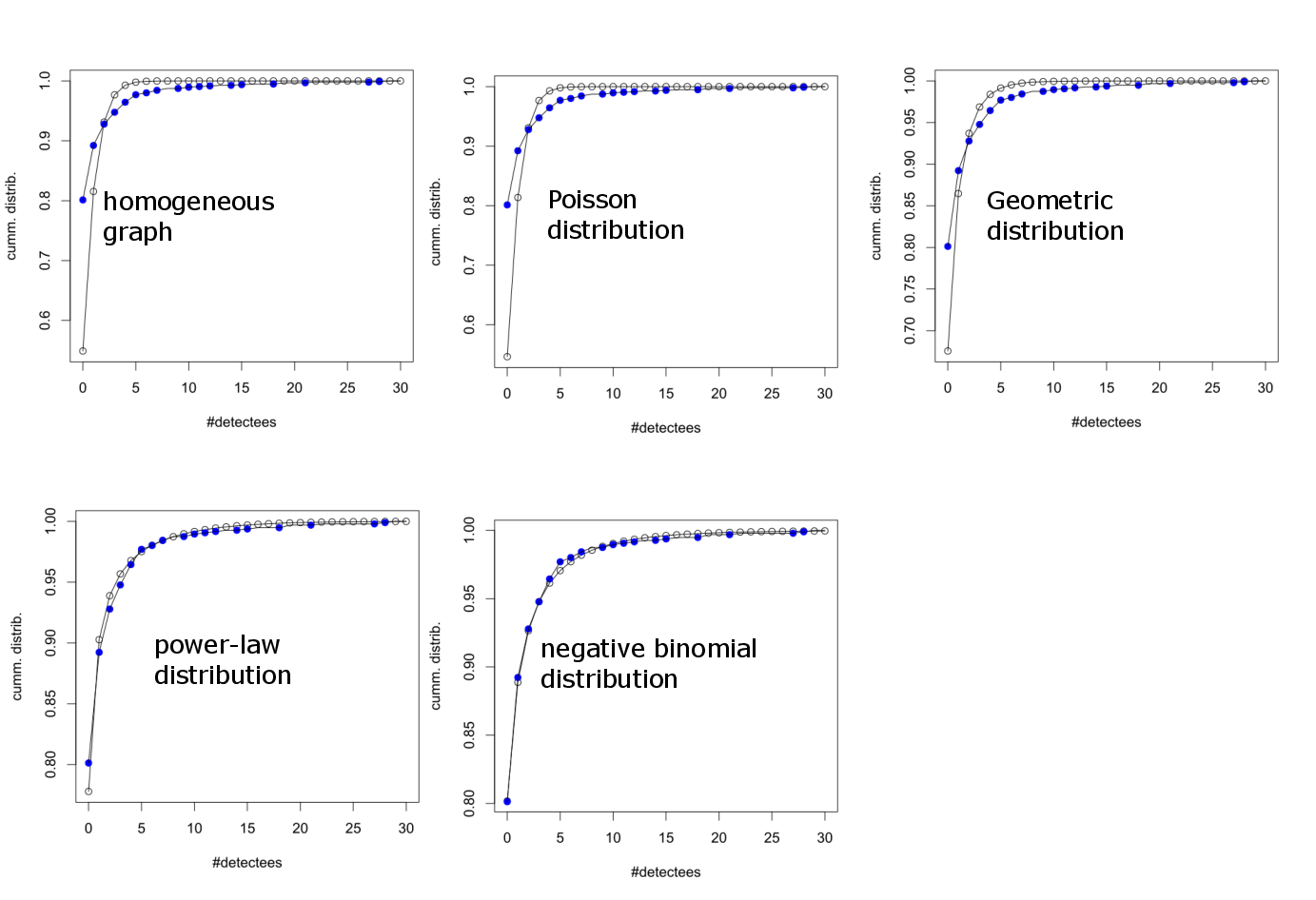}
    \caption{
    Cumulative empirical distribution (blue bullets) compared with the cumulative theoretical distribution (black circles). The top row, from left to right, represents distributions for a homogeneous graph, a Poisson graph, and a geometric distribution, respectively. The bottom row, from left to right, represents distributions for a power-law distribution and a negative binomial distribution, respectively. Our proposed theory is indicated by the black circles.
    }
    \label{fig:cd}
\end{figure}

We first draw a contour plot (fig. \ref{fig:contour_plot}) indicating the point estimates ($p$ and $\mathbb{E}[K]$) and also draw the cumulative empirical distribution vs.\ the cumulative theoretical distribution (fig. \ref{fig:cd}), that is, on the $x$-axis, we plot the number of infectees, and on the $y$-axis the percentage of index cases which is the number of infected contact persons or less. 

Also, these graphics clearly indicate that the power-law and the negative binomial distribution 
yield the best fit, where the negative binomial distribution is superior to the power-law. 

\section{Sensitivity analysis}
\begin{figure}[h!]
    \centering
    \includegraphics[width=1.0\textwidth]{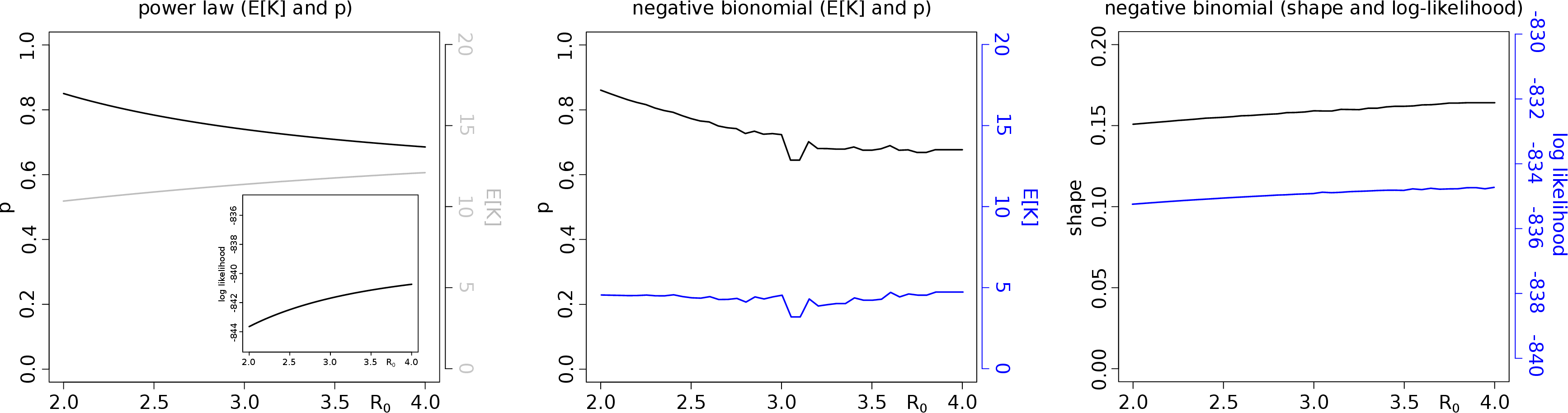}
    \caption{
    Sensitivity analysis illustrating changes in the point estimates ($p$ and $E(K)$), and log-likelihood for different choices of $R_{0}$. Left: Power law (with the log-likelihood as inlay), middle and right: negative binomial distribution.}
    \label{fig:sensi}
\end{figure}

We carried out a sensitivity analysis on how the estimated parameters depend on $R_{0}$. As illustrated in Figure \ref{fig:sensi}, the sensitivity analysis revealed that the estimated parameters in the studied models are not highly sensitive to the choice of $R_{0}$, at least for the power-law and negative distribution (the only degree distributions which meet the data satisfyingly). This observation is particularly noteworthy, as it underscores the robustness of these models in providing reliable estimates of epidemiological parameters, even when the initial assumptions about $R_{0}$ may vary. This characteristic is crucial in the context of real-world epidemiological studies, where the precise value of $R_{0}$ is often uncertain due to factors such as heterogeneous populations, changing contact patterns, and varying degrees of intervention measures. 
 
\section{Discussion}

In this paper, we present a graph-based method for estimating parameters in infectious disease models, offering valuable insights into the efficiency of contact tracing programs and some information about local contact structures, and their implications for the spread of infectious diseases. By comparing various degree distributions and assessing their suitability for modelling disease spread, our analysis contributes to the ongoing development of improved parameter estimation techniques in graph-based models. 
Our findings complement and expand upon the work of~\citet{muller2007estimating} who estimated tracing probability in homogeneous populations by a maximum likelihood estimator and applied it to tuberculosis and chlamydia contact tracing data.
In contrast to focusing on homogeneous random mixing populations, our work explores contact graphs in the form of trees which enables us to capture the unique branching structure of infection transmission chains. This approach provides a more realistic representation of contact patterns at a microscopic level, allowing for a better understanding of the dynamics of infectious diseases and the effectiveness of contact tracing strategies. 
\par\medskip
Our comprehensive analysis of the COVID-19 contact tracing data from Karnataka, India, reveals that both scale-free network models and negative binomial distribution models offer a good fit for the data. The negative binomial distribution emerges as the most fitting model for the data, aligning with previous epidemiological research that has identified this distribution as a suitable assumption for the number of contacts relevant for the transmission of respiratory infections~\cite{mossong2008social}. Furthermore, the observed over-dispersion in the number of secondary cases caused by individual index cases is accurately captured by the negative binomial distribution. This distribution is suitable for data where the variance exceeds the mean, reflective of scenarios where a small proportion of index cases are responsible for a disproportionate number of secondary infections. These findings resonate with the work of \citet{gupta2022contact} who had previously reported a clear over-dispersion in the data. Specifically, \citet{gupta2022contact} found that amongst 956 confirmed index cases, just 8.7\% of cases, who had 14.4\% of contacts, were responsible for 80\% of all secondary cases. The power-law distribution also offers a reasonable fit, highlighting the potential relevance of scale-free networks in modelling infectious disease dynamics. In line with the principles of scale-free networks, our model highlights the role of a relatively small number of "super-spreader" individuals, who have a significantly larger number of contacts and thus a higher likelihood of transmitting the infection to a larger pool of people. This also validates the findings of \citet{gupta2022contact} who suggested that super-spreaders may have played a more dominant role in the covid-19 transmission in Karnataka, India.
\par\medskip
Both the scale-free and negative binomial models allow for a thick or heavy tail, which is created by super-spreader events and it is known in the case of airborne infection such as covid-19 that super-spreader events and the over-dispersion of secondary cases have a significant impact on the effectiveness of contact tracing and surveillance schemes ~\cite{lloyd2005superspreading, gupta2022contact}. Regarding contact tracing, our models consistently indicate a high probability of successful tracing. These figures indeed raise questions regarding their realism and implications for epidemic dynamics. High probabilities for successful contact tracing suggest efficient public health measures in place, as well as the robustness of contact networks to facilitate tracing. However, a high tracing success probability which is not reflective of the high frequency of zero and few traced cases in the reported data might also raise concerns. In particular, regarding the choice of a tree-network where it is assumed that all infections present in the data are part of a single transmission tree. In real-world contact networks where individuals may be part of multiple overlapping transmission chains,  contacts from outside the assumed tree are counted as part of it, potentially resulting in a higher number of traced contacts than what might occur in a non-tree-like network.
\par\medskip
However, the high frequency of zero and few traced cases may have reflected the predominance of cases with a younger age-since-infection during the first wave due to certain preventive measures. Karnataka's contact tracing system was one of India’s most effective, at least, during the early epidemic \cite{team2020laboratory}. Considering the large proportion of close positive asymptomatic contacts at the time of testing, the low numbers present in the data could be indicative of effective pre-testing and preventive strategies at play. For instance, index cases identified and isolated quickly due to effective social distancing and lockdown measures, may result in fewer or no infectious secondary contacts, consequently leading to fewer or no traced cases. This may not necessarily reflect the efficacy of the contact tracing process itself but rather successful containment and prevention efforts that halted the spread from those index cases. This lack of distinction in the overall pre- and post-control strategies in our theory may have overestimated the tracing probability.
\par\medskip
Furthermore, other epidemiological metric such as the reproduction number $R_{0}$ is critical to estimating the number of contacts to trace \cite{linka2020reproduction}. For instance, \citet{hellewell2020feasibility}, who used simulations to study the feasibility of controlling COVID-19 outbreaks through the isolation of cases and contacts, found that to control 90\% of outbreaks with a reproduction number of 2.5, 80\% of contacts needed to be traced and isolated. Their research highlights the subtle role of reproduction number in contact tracing success, revealing that the probability of control increases at all levels of contact tracing when the reproduction number is reduced. Moreover, the \citet{hellewell2020feasibility} study emphasizes the significant impact of the number of initial cases on the likelihood of achieving control. Such insights highlight the complexities of contact tracing and its dependencies on various epidemiological and social factors. Additionally, while the high probabilities of contact tracing success implied by our model may raise questions about their realism, these figures are not unfounded. Our estimator also considers the interplay between the probability of tracing a contact once an index case is identified ($p$) and the probability of an infected individual being detected as an index case ($p_{obs}$). The latter has to be inherently lower, especially within close-knit contact networks such as family units. 
\\
In such networks, it is improbable for all contacts to become index cases; instead, tracing often occurs through one or a few known cases, underscoring why $p_{obs}$ has to be substantially small. Given the nature of our methodology and the reported data, which concentrates on tracing only immediate contacts (one-step tracing), a higher $p$ is plausible while maintaining lower $p_{obs}$. This is reflective of an efficient tracing system where immediate contacts are quickly identified, but not all are independently confirmed as index cases due to the close connection and simultaneous discovery through a single or limited number of initial cases. Despite $p$ exceeding the ideal range for our approximations, the considerably smaller probability of any particular infected person becoming an index case ensures that our model's overall estimations are well suited. It is also crucial to understand that these probabilities, while informative, also underline the inherent complexities in predicting real-world outcomes. Variations in regional practices, public response, healthcare infrastructure, and other socio-cultural factors play a significant role in the success of contact tracing endeavours. Therefore, while our model provides an essential tool for estimation, the results should be interpreted in conjunction with the broader epidemiological context and in light of other research findings for a holistic understanding.
\par\medskip
Our modelling study and findings highlight the importance of selecting appropriate models for estimating tracing probabilities and local contact structures in real-world scenarios. These estimates demonstrate the effectiveness of our graph-based method in capturing key epidemic parameters within heterogeneous and age-structured contact networks. The estimated degree ranges for the negative binomial and power-law distributions fall within plausible ranges found in the literature, supporting the validity of our approach in comparison to other studies that have examined contact tracing data~\cite{mossong2008social, soetens2018real, fyles2021using}. Furthermore, sensitivity analysis on the estimated parameters w.r.t.\ $R_{0}$ for the power-law and negative binomial distribution shows limited sensitivity to the choice of $R_{0}$. This provides valuable insights into estimating key epidemiological and intervention parameters with greater confidence, enabling more effective public health strategies and interventions. All in all, our research contributes to the ongoing development of improved parameter estimation techniques in graph-based models for infectious disease dynamics. By utilizing a graph-based approach and building upon the methods of previous studies, we have demonstrated the value of incorporating contact graph structures, such as trees, for a more accurate representation of contact patterns and infectious disease dynamics. The results of our analysis highlight the need to consider heterogeneity in individual-level contact networks when designing and evaluating contact tracing strategies.
\par\medskip
Future research in this area could explore the incorporation of additional data sources and model refinements. A key area of refinement could be incorporating temporal dynamics or considering other types of contact graphs that better represent real-world contact patterns. The use of tree-like graphs in the current study, while mathematically convenient, is a simplification of reality. Contact patterns, particularly within clusters such as households or other social groups, often exhibit significant clustering and interconnectedness that a tree structure may fail to accurately capture. Our simulation study further highlights these concerns. For instance, the influence of triangles on the fraction of outside infections and the difference in the estimator's performance between Poisson and fixed degree models offer key insights. Such findings suggest a well-connected core group in certain models that leads to more outside infections, challenging our estimator's accuracy. Integrating well-suited graph structures with modelling techniques such as agent-based models or compartmental models, could provide a more comprehensive understanding of infectious disease dynamics and inform the design of more effective public health interventions.

\par\bigskip
{\bf Author Contributions}

{\it Augustine Okolie (A.O.) was the primary investigator and lead author of this manuscript. A.O. conceptualized the study design, performed the data analysis, interpreted the results, and drafted the manuscript. Johannes Müller (J.M.) provided extensive guidance and supervision throughout the development of the research, and contributed to refining the study's methodological approach, data analysis and code implementation. J.M. also provided critical revisions to the manuscript, adding intellectual content and ensuring the integrity of the work. Mirjam Kretzschmar (M.K.) provided substantial inputs to the research including refining the study's methodological approach, analysis and interpretation of data, critical review and commentary on the study's findings. All authors discussed the results and implications and commented on the manuscript at all stages.}

\par\bigskip
{\bf Acknowledgements}

{\it Portions of the work presented in this paper have previously been published in the PhD thesis authored by Augustine \citet{okolie2022contact} and supervised by Johannes Müller. The thesis was submitted to the Technical University of Munich and is accessible online at \url{https://mediatum.ub.tum.de/1661774}. This current paper extends and enhances the findings discussed in the thesis, particularly those related to the likelihood estimator models and their application to both simulated data and empirical COVID-19 contact tracing data, as presented in chapters 3 to 5. The research was supported by a grant from the German Academic Exchange Service (DAAD) (AO), and by the Deutsche Forschungsgemeinschaft (DFG) through the TUM International Graduate School of Science and Engineering (IGSSE), GSC 81, as part of the project GENOMIE\_QADOP (JM). M.K. acknowledges the support from the Horizon 2020 research and innovation funding programme (Grant no. 101003480 (CORESMA)).}

\par\bigskip
{\bf Data Accessibility, Ethics}

{\it The analysis presented in this paper utilizes a publicly available dataset on COVID-19 contact tracing in Karnataka, India, originally reported by \citet{gupta2022contact}. The comprehensive description of the dataset, including the source of data, methodology for data compilation, and ethical considerations, has been documented extensively by \citet{gupta2022contact}. We have adhered to all ethical guidelines as laid out by the original study, which includes the appropriate ethical approvals for data collection and usage. Our research did not involve direct data collection from human subjects, and all analysis were conducted on anonymized, aggregated data without any personal identifiers to maintain individual confidentiality. The original dataset was collected in compliance with ethical standards of research, and our use of this dataset for secondary analysis is consistent with those standards.}

\par\bigskip
{\bf Declaration of Competing Interest}
\par\medskip
{None}

\par\bigskip
{\bf Supplementary material}
\par\medskip
Supplementary material associated with this article can be found in the published version at \url{https://doi.org/10.1098/rsif.2023.0409}.
	
\clearpage
\bibliography{references}

\begin{thebibliography}{41}
\providecommand{\natexlab}[1]{#1}
\providecommand{\url}[1]{\texttt{#1}}
\expandafter\ifx\csname urlstyle\endcsname\relax
  \providecommand{\doi}[1]{doi: #1}\else
  \providecommand{\doi}{doi: \begingroup \urlstyle{rm}\Url}\fi

\bibitem[Khan et~al.(2020)]{khan2020parameter}
Muhammad~Altaf Khan et~al.
\newblock Parameter estimation and fractional derivatives of dengue
  transmission model.
\newblock \emph{AIMS Mathematics}, 5\penalty0 (3):\penalty0 2758--2779, 2020.

\bibitem[Little et~al.(2012)Little, D'Agostino, Cohen, Dickersin, Emerson,
  Farrar, Frangakis, Hogan, Molenberghs, Murphy, et~al.]{little2012prevention}
Roderick~J Little, Ralph D'Agostino, Michael~L Cohen, Kay Dickersin, Scott~S
  Emerson, John~T Farrar, Constantine Frangakis, Joseph~W Hogan, Geert
  Molenberghs, Susan~A Murphy, et~al.
\newblock The prevention and treatment of missing data in clinical trials.
\newblock \emph{New England Journal of Medicine}, 367\penalty0 (14):\penalty0
  1355--1360, 2012.

\bibitem[Cant{\'o} et~al.(2009)Cant{\'o}, Coll, and
  S{\'a}nchez]{canto2009structural}
Bego{\~n}a Cant{\'o}, Carmen Coll, and Elena S{\'a}nchez.
\newblock Structural identifiability of a model of dialysis.
\newblock \emph{Mathematical and computer modelling}, 50\penalty0
  (5-6):\penalty0 733--737, 2009.

\bibitem[Cant{\'o} et~al.(2011)Cant{\'o}, Coll, and
  S{\'a}nchez]{canto2011identifiability}
Bego{\~n}a Cant{\'o}, Carmen Coll, and Elena S{\'a}nchez.
\newblock Identifiability for a class of discretized linear partial
  differential algebraic equations.
\newblock \emph{Mathematical Problems in Engineering}, 2011, 2011.

\bibitem[Craciun and Pantea(2008)]{craciun2008identifiability}
Gheorghe Craciun and Casian Pantea.
\newblock Identifiability of chemical reaction networks.
\newblock \emph{Journal of Mathematical Chemistry}, 44\penalty0 (1):\penalty0
  244--259, 2008.

\bibitem[Blum and Tran(2010)]{blum2010hiv}
Michael~GB Blum and Viet~Chi Tran.
\newblock Hiv with contact tracing: a case study in approximate bayesian
  computation.
\newblock \emph{Biostatistics}, 11\penalty0 (4):\penalty0 644--660, 2010.

\bibitem[O'Neill et~al.(1997)O'Neill, Roberts, and
  of~Mathematics;]{o1997bayesian}
PD~O'Neill, GO~Roberts, and Bradford Univ.(United Kingdom).~Dept.
  of~Mathematics;.
\newblock \emph{Bayesian inference for partially observed stochastic
  epidemics}.
\newblock University of Bradford. School of Mathematical Sciences, 1997.

\bibitem[G{\"o}tz et~al.(2017)G{\"o}tz, Altmeier, Bock, Rockenfeller, Wijaya,
  et~al.]{gotz2017modeling}
Thomas G{\"o}tz, Nicole Altmeier, Wolfgang Bock, Robert Rockenfeller,
  Karunia~Putra Wijaya, et~al.
\newblock Modeling dengue data from semarang, indonesia.
\newblock \emph{Ecological complexity}, 30:\penalty0 57--62, 2017.

\bibitem[Agusto and Khan(2018)]{agusto2018optimal}
FB~Agusto and MA~Khan.
\newblock Optimal control strategies for dengue transmission in pakistan.
\newblock \emph{Mathematical biosciences}, 305:\penalty0 102--121, 2018.

\bibitem[Stollenwerk et~al.(2012)Stollenwerk, Aguiar, Ballesteros, Boto, Kooi,
  and Mateus]{Stollenwerk2012}
N.~Stollenwerk, M.~Aguiar, S.~Ballesteros, J.~Boto, B.~Kooi, and L.~Mateus.
\newblock Dynamic noise, chaos and parameter estimation in population biology.
\newblock \emph{Interface Focus}, 2\penalty0 (2):\penalty0 156--169, feb 2012.
\newblock \doi{10.1098/rsfs.2011.0103}.

\bibitem[Manou-Abi et~al.(2022)Manou-Abi, Slaoui, and
  Balicchi]{manou2022estimation}
Solym~M Manou-Abi, Yousri Slaoui, and Julien Balicchi.
\newblock Estimation of some epidemiological parameters with the covid-19 data
  of mayotte.
\newblock \emph{Frontiers in Applied Mathematics and Statistics ()}, page~67,
  2022.

\bibitem[M{\"u}ller and H{\"o}sel(2007)]{muller2007estimating}
Johannes M{\"u}ller and Volker H{\"o}sel.
\newblock Estimating the tracing probability from contact history at the onset
  of an epidemic.
\newblock \emph{Mathematical Population Studies}, 14\penalty0 (4):\penalty0
  211--236, 2007.

\bibitem[M{\"u}ller et~al.(2000)M{\"u}ller, Kretzschmar, and
  Dietz]{muller2000contact}
Johannes M{\"u}ller, Mirjam Kretzschmar, and Klaus Dietz.
\newblock Contact tracing in stochastic and deterministic epidemic models.
\newblock \emph{Mathematical biosciences}, 164\penalty0 (1):\penalty0 39--64,
  2000.

\bibitem[Dyson et~al.(2018)Dyson, Marks, Crook, Sokana, Solomon, Bishop, Mabey,
  and Hollingsworth]{dyson2018targeted}
Louise Dyson, Michael Marks, Oliver~M Crook, Oliver Sokana, Anthony~W Solomon,
  Alex Bishop, David~CW Mabey, and T~D{\'e}irdre Hollingsworth.
\newblock Targeted treatment of yaws with household contact tracing: how much
  do we miss?
\newblock \emph{American journal of epidemiology}, 187\penalty0 (4):\penalty0
  837--844, 2018.

\bibitem[Tanaka et~al.(2020)Tanaka, Yamaguchi, and Sakamoto]{Tanaka2020}
Takuma Tanaka, Takayuki Yamaguchi, and Yohei Sakamoto.
\newblock Estimation of the percentages of undiagnosed patients of the novel
  coronavirus ({SARS}-{CoV}-2) infection in hokkaido, japan by using
  birth-death process with recursive full tracing.
\newblock \emph{{PLOS} {ONE}}, 15:\penalty0 e0241170, 2020.
\newblock \doi{10.1371/journal.pone.0241170}.

\bibitem[Okolie and M{\"u}ller(2020)]{okolie2020exact}
Augustine Okolie and Johannes M{\"u}ller.
\newblock Exact and approximate formulas for contact tracing on random trees.
\newblock \emph{Mathematical biosciences}, 321:\penalty0 108320, 2020.

\bibitem[Bertacchini et~al.(2020)Bertacchini, Bilotta, and
  Pantano]{bertacchini2020temporal}
Francesca Bertacchini, Eleonora Bilotta, and Pietro~S Pantano.
\newblock On the temporal spreading of the sars-cov-2.
\newblock \emph{PLoS One}, 15\penalty0 (10):\penalty0 e0240777, 2020.

\bibitem[Chondros et~al.(2022)Chondros, Nikolopoulos, and
  Polenakis]{chondros2022integrated}
Christos Chondros, Stavros~D Nikolopoulos, and Iosif Polenakis.
\newblock An integrated simulation framework for the prevention and mitigation
  of pandemics caused by airborne pathogens.
\newblock \emph{Network Modeling Analysis in Health Informatics and
  Bioinformatics}, 11\penalty0 (1):\penalty0 42, 2022.

\bibitem[Cuevas-Maraver et~al.(2021)Cuevas-Maraver, Kevrekidis, Chen,
  Kevrekidis, Villalobos-Daniel, Rapti, and Drossinos]{cuevas2021lockdown}
Jes{\'u}s Cuevas-Maraver, Panayotis~G Kevrekidis, Qian-Yong Chen, George~A
  Kevrekidis, V{\'\i}ctor Villalobos-Daniel, Zoi Rapti, and Yannis Drossinos.
\newblock Lockdown measures and their impact on single-and two-age-structured
  epidemic model for the covid-19 outbreak in mexico.
\newblock \emph{Mathematical Biosciences}, 336:\penalty0 108590, 2021.

\bibitem[Kovacevic et~al.(2020)Kovacevic, Stilianakis, and
  Veliov]{kovacevic2020distributed}
Raimund Kovacevic, Nikolaos~I Stilianakis, and Vladimir~M Veliov.
\newblock A distributed optimal control epidemiological model applied to
  covid-19 pandemic.
\newblock \emph{available as ORCOS Research Report}, 13, 2020.

\bibitem[Modi et~al.(2021)Modi, Umate, Makade, Dubey, and
  Agarwal]{modi2021simulation}
Kanak Modi, Laxmikant Umate, Kiran Makade, Ravi~Shanker Dubey, and Pankaj
  Agarwal.
\newblock Simulation based study for estimation of covid-19 spread in india
  using seir model.
\newblock \emph{Journal of Interdisciplinary Mathematics}, 24\penalty0
  (2):\penalty0 245--258, 2021.

\bibitem[Kevrekidis et~al.(2021)Kevrekidis, Cuevas-Maraver, Drossinos, Rapti,
  and Kevrekidis]{kevrekidis2021reaction}
Panayotis~G Kevrekidis, Jes{\'u}s Cuevas-Maraver, Yannis Drossinos, Zoi Rapti,
  and George~A Kevrekidis.
\newblock Reaction-diffusion spatial modeling of covid-19: Greece and andalusia
  as case examples.
\newblock \emph{Physical Review E}, 104\penalty0 (2):\penalty0 024412, 2021.

\bibitem[Green and Kiss(2010)]{green2010large}
Darren~M Green and Istvan~Z Kiss.
\newblock Large-scale properties of clustered networks: Implications for
  disease dynamics.
\newblock \emph{Journal of biological dynamics}, 4\penalty0 (5):\penalty0
  431--445, 2010.

\bibitem[Keeling(1999)]{keeling1999effects}
Matthew~J Keeling.
\newblock The effects of local spatial structure on epidemiological invasions.
\newblock \emph{Proceedings of the Royal Society of London. Series B:
  Biological Sciences}, 266\penalty0 (1421):\penalty0 859--867, 1999.

\bibitem[Kiss et~al.(2017)Kiss, Miller, Simon, et~al.]{kiss2017mathematics}
Istv{\'a}n~Z Kiss, Joel~C Miller, P{\'e}ter~L Simon, et~al.
\newblock Mathematics of epidemics on networks.
\newblock \emph{Cham: Springer}, 598:\penalty0 31, 2017.

\bibitem[Newman(2002)]{newman2002spread}
Mark~EJ Newman.
\newblock Spread of epidemic disease on networks.
\newblock \emph{Physical review E}, 66\penalty0 (1):\penalty0 016128, 2002.

\bibitem[Karrer and Newman(2010)]{karrer2010message}
Brian Karrer and Mark~EJ Newman.
\newblock Message passing approach for general epidemic models.
\newblock \emph{Physical Review E}, 82\penalty0 (1):\penalty0 016101, 2010.

\bibitem[Feld(1991)]{Feld1991}
Scott~L. Feld.
\newblock Why your friends have more friends than you do.
\newblock \emph{American Journal of Sociology}, 96\penalty0 (6):\penalty0
  1464--1477, may 1991.
\newblock \doi{10.1086/229693}.

\bibitem[Brauer et~al.(2008)Brauer, Van~den Driessche, and
  Wu]{allen2008mathematical}
Fred Brauer, Pauline Van~den Driessche, and Jianhong Wu.
\newblock \emph{Mathematical epidemiology}, volume 1945.
\newblock Springer, 2008.

\bibitem[Christakis and Fowler(2010)]{Christakis2010}
Nicholas~A. Christakis and James~H. Fowler.
\newblock Social network sensors for early detection of contagious outbreaks.
\newblock \emph{{PLoS} {ONE}}, 5\penalty0 (9):\penalty0 e12948, sep 2010.
\newblock \doi{10.1371/journal.pone.0012948}.

\bibitem[Lloyd-Smith et~al.(2005)Lloyd-Smith, Schreiber, Kopp, and
  Getz]{lloyd2005superspreading}
James~O Lloyd-Smith, Sebastian~J Schreiber, P~Ekkehard Kopp, and Wayne~M Getz.
\newblock Superspreading and the effect of individual variation on disease
  emergence.
\newblock \emph{Nature}, 438\penalty0 (7066):\penalty0 355--359, 2005.

\bibitem[Gupta et~al.(2022)Gupta, Parameswaran, Sra, Mohanta, Patel, Gupta,
  Bansal, Jain, Mazumder, Arora, et~al.]{gupta2022contact}
Mohak Gupta, Giridara~G Parameswaran, Manraj~S Sra, Rishika Mohanta, Devarsh
  Patel, Amulya Gupta, Bhavik Bansal, Vardhmaan Jain, Archisman Mazumder, Mehak
  Arora, et~al.
\newblock Contact tracing of covid-19 in karnataka, india: Superspreading and
  determinants of infectiousness and symptomatic infection.
\newblock \emph{Plos one}, 17\penalty0 (7):\penalty0 e0270789, 2022.

\bibitem[Bell et~al.(1994)Bell, Goldoft, Griffin, Davis, Gordon, Tarr,
  Bartleson, Lewis, Barrett, Wells, et~al.]{bell1994multistate}
Beth~P Bell, Marcia Goldoft, Patricia~M Griffin, Margaret~A Davis, Diane~C
  Gordon, Phillip~I Tarr, Charles~A Bartleson, Jay~H Lewis, Timothy~J Barrett,
  Joy~G Wells, et~al.
\newblock A multistate outbreak of escherichia coli o157: h7—associated
  bloody diarrhea and hemolytic uremic syndrome from hamburgers: the washington
  experience.
\newblock \emph{Jama}, 272\penalty0 (17):\penalty0 1349--1353, 1994.

\bibitem[Locatelli et~al.(2021)Locatelli, Tr{\"a}chsel, and
  Rousson]{locatelli2021estimating}
Isabella Locatelli, Bastien Tr{\"a}chsel, and Valentin Rousson.
\newblock Estimating the basic reproduction number for covid-19 in western
  europe.
\newblock \emph{Plos one}, 16\penalty0 (3):\penalty0 e0248731, 2021.

\bibitem[Mossong et~al.(2008)Mossong, Hens, Jit, Beutels, Auranen, Mikolajczyk,
  Massari, Salmaso, Tomba, Wallinga, et~al.]{mossong2008social}
Jo{\"e}l Mossong, Niel Hens, Mark Jit, Philippe Beutels, Kari Auranen, Rafael
  Mikolajczyk, Marco Massari, Stefania Salmaso, Gianpaolo~Scalia Tomba, Jacco
  Wallinga, et~al.
\newblock Social contacts and mixing patterns relevant to the spread of
  infectious diseases.
\newblock \emph{PLoS medicine}, 5\penalty0 (3):\penalty0 e74, 2008.

\bibitem[Team et~al.(2020)Team, Team, Abraham, Aggarwal, Babu, Barani,
  Bhargava, Bhatnagar, Dhama, Gangakhedkar, et~al.]{team2020laboratory}
Data~Management Team, COVID~Laboratory Team, Priya Abraham, Neeraj Aggarwal,
  Giridhara~R Babu, Suganya Barani, Balram Bhargava, Tarun Bhatnagar,
  Ajay~Singh Dhama, Raman~R Gangakhedkar, et~al.
\newblock Laboratory surveillance for sars-cov-2 in india: Performance of
  testing \& descriptive epidemiology of detected covid-19, january 22-april
  30, 2020.
\newblock \emph{The Indian journal of medical research}, 151\penalty0
  (5):\penalty0 424, 2020.

\bibitem[Linka et~al.(2020)Linka, Peirlinck, and Kuhl]{linka2020reproduction}
Kevin Linka, Mathias Peirlinck, and Ellen Kuhl.
\newblock The reproduction number of covid-19 and its correlation with public
  health interventions.
\newblock \emph{Computational mechanics}, 66:\penalty0 1035--1050, 2020.

\bibitem[Hellewell et~al.(2020)Hellewell, Abbott, Gimma, Bosse, Jarvis,
  Russell, Munday, Kucharski, Edmunds, Sun, et~al.]{hellewell2020feasibility}
Joel Hellewell, Sam Abbott, Amy Gimma, Nikos~I Bosse, Christopher~I Jarvis,
  Timothy~W Russell, James~D Munday, Adam~J Kucharski, W~John Edmunds, Fiona
  Sun, et~al.
\newblock Feasibility of controlling covid-19 outbreaks by isolation of cases
  and contacts.
\newblock \emph{The Lancet Global Health}, 8\penalty0 (4):\penalty0 e488--e496,
  2020.

\bibitem[Soetens et~al.(2018)Soetens, Klinkenberg, Swaan, Hahn{\'e}, and
  Wallinga]{soetens2018real}
Loes Soetens, Don Klinkenberg, Corien Swaan, Susan Hahn{\'e}, and Jacco
  Wallinga.
\newblock Real-time estimation of epidemiologic parameters from contact tracing
  data during an emerging infectious disease outbreak.
\newblock \emph{Epidemiology}, 29\penalty0 (2):\penalty0 230--236, 2018.

\bibitem[Fyles et~al.(2021)Fyles, Fearon, Overton, of~Manchester COVID-19
  Modelling~Group, Wingfield, Medley, Hall, Pellis, and House]{fyles2021using}
Martyn Fyles, Elizabeth Fearon, Christopher Overton, University of~Manchester
  COVID-19 Modelling~Group, Tom Wingfield, Graham~F Medley, Ian Hall, Lorenzo
  Pellis, and Thomas House.
\newblock Using a household-structured branching process to analyse contact
  tracing in the sars-cov-2 pandemic.
\newblock \emph{Philosophical Transactions of the Royal Society B},
  376\penalty0 (1829):\penalty0 20200267, 2021.

\bibitem[Okolie(2022)]{okolie2022contact}
Augustine~Okebunor Okolie.
\newblock \emph{Contact Tracing on Stochastic Graphs}.
\newblock PhD thesis, Universit{\"a}tsbibliothek der TU M{\"u}nchen, 2022.

\end{thebibliography}
\bibliographystyle{unsrtnat}

\newpage
\begin{appendix}
\section{Full graph/random mixing}\label{randMixAppend}
In the case of a full graph, we start with a fixed degree $K=N-1$, where $N$ is the population size, and take the limit $N\rightarrow\infty$ and $\beta\rightarrow 0$ such that $R_0=\frac{(N-1)\,\beta}{\mu+\sigma}$ is constant. 
Then,  
$$ (N-1)\, \hat p(a) 
= R_0\, \frac{\alpha+\sigma}{\alpha+\sigma-\beta}\, \bigg(e^{-\beta a}-e^{-(\alpha+\sigma) a}\bigg)
\rightarrow 
R_0\, (1-e^{-(\alpha+\sigma)\,a}). 
$$
Therefore the Binomial distribution approximates a Poisson distribution, and eqn~\ref{T_result} becomes in the limit $N\rightarrow\infty$ 
$$ P(T=i) =
 \int_{0}^{\infty} \,\, \mbox{dpois}(i, R_0\, (1-e^{-(\alpha+\sigma)\,a})\,\,)\,\,\,\,R_0\, (\alpha+\sigma)\, e^{-R_0\, (\alpha+\sigma) a}\, da
$$
where $\mbox{dpois}(i,\mu) = \mu^i e^{-\mu}/i!$ is the probability function 
for the Poisson distribution.

\section{Optimization Process}\label{Optim}

Given a function $f(x)$, the aim of optimization is to find an $x$ that either maximizes or minimizes $f(x)$. This process involves the computation of gradients and the Hessian matrix.

The gradient of a function is a vector that points in the direction of the greatest increase of that function. It is calculated as the vector of the first derivatives of the function with respect to each variable. The gradient of a function $f(x)$, where $x = [x_1, x_2, \ldots, x_n]$ is:

\begin{equation}
\nabla f(x) = \left[\frac{\partial f}{\partial x_1}, \frac{\partial f}{\partial x_2}, \ldots, \frac{\partial f}{\partial x_n}\right]
\label{eqn:grad}
\end{equation}

To ensure that the solution found is a local maximum and not a local minimum or a saddle point, the Hessian matrix is used. The Hessian matrix is the square matrix of second-order partial derivatives of the function. Each element in the Hessian matrix is the second derivative of the function with respect to different variables. The Hessian matrix for the function $f(x)$ is:

\begin{equation}
H(f(x)) = 
\begin{bmatrix}
\frac{\partial^2 f}{\partial x_1^2} & \frac{\partial^2 f}{\partial x_1\partial x_2} & \dots & \frac{\partial^2 f}{\partial x_1\partial x_n} \\
\frac{\partial^2 f}{\partial x_2\partial x_1} & \frac{\partial^2 f}{\partial x_2^2} & \dots & \frac{\partial^2 f}{\partial x_2\partial x_n} \\
\vdots & \vdots & \ddots & \vdots \\
\frac{\partial^2 f}{\partial x_n\partial x_1} & \frac{\partial^2 f}{\partial x_n\partial x_2} & \dots & \frac{\partial^2 f}{\partial x_n^2} \\
\end{bmatrix}
\end{equation}

If the Hessian is positive definite (all eigenvalues are positive) at a point, then the function attains a local minimum at that point. If the Hessian is negative definite (all eigenvalues are negative), then the function attains a local maximum.

In the context of maximizing a likelihood function, we often convert the problem into a minimization problem by taking the negative of the likelihood function. This is due to the fact that many optimization algorithms are developed for minimization problems. The negative log-likelihood function becomes:

\begin{equation}
- \mathcal{LL} \left(\vec{\mu} \mid i_{\ell}, \ell=1, ..., n\right)
\end{equation}

The goal now is to minimize this negative log-likelihood function, and the optimization problem becomes:

\begin{equation}
\vec{\mu}^{*} = \underset{\vec{\mu}}{\text{argmin}} \left\{- \mathcal{LL} \left(\vec{\mu} \mid i_{\ell}, \ell=1, ..., n\right)\right\}
\end{equation}

The same principles of gradients and Hessians apply to this minimization problem. The gradient of the negative log-likelihood function should point in the direction of the greatest decrease of the function. The Hessian, on the other hand, should be negative definite at the point of minimum.

In the case of the log-likelihood function, the optimization problem can be solved using iterative methods such as Newton's method or Quasi-Newton methods, which make use of both the gradient and the Hessian of the function to find the minimum.
\end{appendix}

\end{document}